\documentclass[11 pt, reqno]{article}

\usepackage{amsmath,amsfonts,amssymb,amsthm}
\usepackage[colorlinks=true,hyperindex=true]{hyperref}
\usepackage{cancel,bbm}
\usepackage{mathabx} 
\usepackage{comment}
\usepackage{enumitem}
\usepackage{cite}
\usepackage{showlabels}
\usepackage{graphicx}

\usepackage{chngcntr}
\counterwithin*{equation}{section}

\usepackage[left=3cm, right=3cm, top=3 cm, bottom= 3 cm]{geometry}
\usepackage{color}
\usepackage{fancyhdr}
\usepackage{latexsym}

\usepackage{bm}

\newtheorem{ccounter}{ccounter}[section]
\newtheorem{thm}[ccounter]{Theorem}
\newtheorem{lem}[ccounter]{Lemma}
\newtheorem{cor}[ccounter]{Corollary}
\newtheorem{defn}[ccounter]{Definition}
\newtheorem{prop}[ccounter]{Proposition}
\newtheorem{ass}[ccounter]{Assumption}
\newtheorem{ex}[ccounter]{Example}

\def\bet{\begin{thm}}
\def\eet{\end{thm}}
\def\bel{\begin{lem}}
\def\eel{\end{lem}}
\def\bas{\begin{ass}}
\def\eas{\end{ass}}
\def\bec{\begin{cor}}
\def\eec{\end{cor}}
\def\bed{\begin{defn}}
\def\eed{\end{defn}}
\def\bep{\begin{prop}}
\def\eep{\end{prop}}
\def\beq{\begin{equation}}
\def\eeq{\end{equation}}
\def\proof{\noindent {\bf Proof.}\ \ }
\def\bea{\begin{equation*}}
\def\eea{\end{equation*}}

\def\bex{\begin{ex}}
\def\eex{\end{ex}}

\def\rr{\mathbb{R}}
\def\zz{\mathbb{Z}}

\def\1{\boldsymbol{1}}

\def\e{\mathrm{e}}
\def\i{\mathrm{i}}
\def\del{\partial}
\def\d{\mathrm{d}}
\def\eps{\varepsilon}
\renewcommand\leq\varleq
\renewcommand\geq\vargeq
\def\ee{\mathrm{E}}

\def\O{\mathcal{O}}

\def\ee{\mathbb{E}}

\def\pp{\mathbb{P}}

\def\I{\mathcal{I}}
\def\mfa{\mathfrak{m}}
\def\A{\mathcal{A}}

\def\mfa{\mathfrak{a}}

\def\B{\mathcal{B}}

\DeclareMathOperator{\ad}{ad}
\DeclareMathOperator{\TF}{TF}
\DeclareMathOperator{\Cov}{Cov}

\begin{document}

\begin{table}
\centering

\begin{tabular}{c}

\multicolumn{1}{c}{\parbox{12cm}{\begin{center}\Large{\bf Tail bounds for the O'Connell-Yor polymer }\end{center}}}\\
\\
\end{tabular}
\begin{tabular}{ c c c  }

 Benjamin Landon &  &  Philippe Sosoe
 \\
 & & \\  
  \footnotesize{University of Toronto} & \footnotesize{} & \footnotesize{Cornell University}  \\
  \footnotesize{Department of Mathematics}  & \footnotesize{} & \footnotesize{Department of Mathematics}  \\
 \footnotesize{\texttt{blandon@math.toronto.edu}} & \footnotesize{\texttt{}} &\footnotesize{\texttt{psosoe@math.cornell.edu}} \\
  & & \\
\end{tabular}
\\
\begin{tabular}{c}
\multicolumn{1}{c}{\today}\\
\\
\end{tabular}

\begin{tabular}{p{15 cm}}
\small{{\bf Abstract:} We derive upper and lower bounds for the upper and lower tails of the O'Connell-Yor polymer of the correct order of magnitude via probabilistic and geometric techniques in the moderate deviations regime.  The inputs of our work are an identity for the generating function of a two-parameter model of Rains and Emrah-Janjigian-Sepp{\"a}l{\"a}inen,  and the geometric techniques of Ganguly-Hegde and Basu-Ganguly-Hammond-Hegde. As an intermediate result we obtain strong tail estimates for the transversal fluctuation of the polymer path from the diagonal. }
\end{tabular}
\end{table}

\tableofcontents

\section{Introduction}
The O'Connell-Yor polymer, also known as the semi-discrete directed polymer, is a fundamental example of a directed polymer in a random environment in 1+1 dimensions, a collection of models that are expected to lie in the Kardar-Parisi-Zhang (KPZ) universality class. It was introduced by O'Connell and Yor \cite{OY} as a positive temperature analog of the much-studied  Brownian Last Passage percolation. These authors studied a stationary version of the model, and showed that it possesses the Burke property, an invariance property found in certain queueing models, which translates to a two-dimensional invariance property in the corresponding polymer and last passage percolation models. Starting with O'Connell \cite{oconnell}, it was later discovered that the model has an even richer underlying integrable structure, a fact which, for the first time in any polymer model,  ultimately enabled the verification of KPZ type asymptotics for the distribution of the normalized free energy, in the breakthrough work of Borodin, Corwin and Ferrari \cite{BCF}. 

In parallel to the development of integrable probability, there has been an increasing interest in geometric and probabilistic methods to analyze the fluctuations of models expected to be in the KPZ class. In the case of the O'Connell-Yor polymer, this line of research was initiated by Sepp\"al\"ainen and Valk\'o \cite{SV}, who first applied the coupling method of Bal{\'a}zs-Cator-Sepp\"al\"ainen \cite{BCS} to this model to obtain cube root scaling for the fluctuations at the level of the variance. See also \cite{MFSV} for the intermediate disorder case, where the variance of the environment (the ``temperature parameter'') is allowed to depend on the system size. In \cite{NS1}, the second author and Noack estimated the higher moments on a near-optimal scale. Although it features Gaussian integration by parts prominently, the method introduced in \cite{NS1} in fact extends to discrete polymer models \cite{NS2}. Most recently, the authors of the current article obtained upper and lower bounds for the upper tail of the stationary O'Connell-Yor polymer and the four integrable discrete polymer models \cite{LS} (see also the thesis of Xie  containing similar results for the discrete models \cite{X}).

Geometric and probabilistic methods have generally not yet been able to provide as detailed information as those of integrable probability. For example, identifying asymptotic distributions without resorting to explicit formulas remains an outstanding challenge. Moreover, current implementations still require some modest integrable inputs like stationarity or the Burke property. However, the gap has been closing. We mention in this context the important recent results of Emrah, Janjigian, and Sepp\"al\"ainen \cite{EJS}, who introduced a methodology for stationary models that allows them to obtain the exact upper tail (including constants in the exponent) for exponential last passage percolation. See also \cite{BGHH,GH,GH2,LS}. The reason these methods are of great interest is that they have the potential to be more robust than integrable methods under perturbations of parameters, including the initial data and, ultimately, the distribution of the underlying environment variables.

In this paper, we deal with a question that has attracted much recent attention, namely the tail behavior of models in the Kardar-Parisi-Zhang universality class. We complement our previous results on the upper tails of the OY polymer with matching upper and lower bounds on the more delicate lower tail. Limiting one-point distributions in the KPZ class, such as the Tracy-Widom and Baik-Rains distribution, exhibit characteristic super-exponential  decay with specific exponents $\frac{3}{2}$ (for the upper tail) and $3$ (for the lower tail). In many cases, KPZ models reproduce this tail behavior, at least in the moderate deviation range, in pre-limiting regimes. Moreover, for some models it is known that the tail exponents remain the same under perturbations of initial conditions and even the form of the model (see \cite{LS}). Tail behavior is a robust characteristic of the KPZ universality class, and the current work is thus a contribution towards a better understanding of this universality.

\subsection{Definition of the model}
The partition function of the O'Connell-Yor polymer is defined by,
\beq
Z_{n, t} := \int_{0 < s_1 < \dots < s_{n-1} < t } \exp \left( \sum_{i=1}^n B_i (s_i ) - B_i (s_{i-1} ) \right) \d s_1 \dots \d s_{n-1},
\eeq
with the convention $s_0 = 0$ and $s_n = t$.  Above, $\{ B_i  \}_{i}$ are a family of standard Brownian motions. Moriarty and O'Connell \cite{MO} calculated the limiting free energy density, for any $t>0$,
\beq
\lim_{n \to \infty} \frac{1}{n} \log Z_{n, tn} = \theta t - \psi_0 ( \theta)
\eeq
where,
\beq
\psi_k ( \theta) :=\frac{ \d^{k+1} }{ \d \theta^{k+1}} \log \left( \int_0^\infty s^{\theta-1} \e^{ -s } \d s \right)
\eeq
are the polygamma functions and $\theta$ is the unique solution to $\psi_1 ( \theta) = t$. In \cite{LS}, we proved the estimate
\beq
\pp\left[ \log Z_{n, tn} > n( \theta t - \psi_0 ( \theta ) ) + s n^{1/3} \right] \leq C \e^{ -c s^{3/2} }, \qquad t = \psi_1 ( \theta )
\eeq
for some constants $c, C>0$ and any $0 < s < c n^{2/3}$.  In the present work, we will complement this upper bound for the upper tail with lower bound of matching order, as well as upper and lower bounds for the lower tail.  This is the content of the following theorem, which summarizes these statements.

\bet
Let $\delta >0$ and assume $\delta n \leq t \leq \delta^{-1} n$. There are $C, c>0$ so that the following hold. Let $\theta$ satisfy $\psi_1 ( \theta ) = t/n$.  We have,
\beq \label{eqn:main-upper-tail}
c \e^{ - C s^{3/2} } \leq \pp\left[ \log Z_{n, t} > ( \theta t - n\psi_0 ( \theta ) ) + s n^{1/3} \right]  \leq C \e^{ - c s^{3/2}}
\eeq
for all $0 < s < c n^{2/3}$.  For $0 < s < c n^{2/3} / \log(n)$ we have,
\beq \label{eqn:main-lower-upper}
\pp\left[ \log Z_{n, t} <  ( \theta t-n \psi_0 ( \theta ) )- s n^{1/3} \right] \leq C \e^{ - cs^3}
\eeq
and for $0 < s < c n^{2/3} / \log(n)^2$ we have,
\beq \label{eqn:main-lower-lower}
\pp\left[ \log Z_{n, t} <  ( \theta t-n \psi_0 ( \theta ) )- s n^{1/3} \right] \geq c \e^{ - C s^3}.
\eeq
\eet
The various estimates of the above theorem are proven in the following sections of the paper. 
The upper bound of \eqref{eqn:main-upper-tail} follows from Proposition \ref{prop:old-tail} (a restatement of results of \cite{LS}), and the lower bound is proven in Section \ref{sec:upper-tail}. 
The estimate \eqref{eqn:main-lower-upper} follows from Proposition \ref{prop:watermelon-bound}. The estimate \eqref{eqn:main-lower-lower} follows from Theorem \ref{thm:lt-lower}. 

In addition to the moderate deviations tail estimates for exponential last passage percolation of Emrah, Janjigian and Sepp{\"a}l{\"a}inen \cite{EJS} mentioned in the introduction, we also mention the related large deviations estimates for the O'Connell-Yor polymer that were proven by Janjigian \cite{J}. This work builds on the approach of \cite{GS} for the log gamma polymer. This corresponds to the regime $s = \O ( n^{2/3})$ where the KPZ tail exponents are no longer expected to arise.

\subsection{Methodology}

In the work \cite{LS} we considered stationary KPZ models (the stationary OY polymer is introduced in Section \ref{sec:models} below) and showed how monotonicity and convexity of the models in the parameters defining the systems and a certain identity involving the moment generating function of a two-parameter version of the model could be combined to yield a short and transparent proof of an upper bound for the lower and upper tails of the form $\e^{ - c s^{3/2}}$ for the stationary models. This identity was first derived by Rains \cite{Rains} in the context of last passage percolation, but was recently re-introduced and used to great effect in the work of Emrah-Janjigian-Sepp{\"a}l{\"a}inen \cite{EJS}. Due to the fact that the two-parameter model stochastically dominates the non-stationary model $\log Z_{n, t}$ considered here, the Rains-EJS identity in fact yields a short proof of the upper bound for the upper tail, as indicated in \cite{LS}. 

The main contributions of this work are then  the remaining estimates, the lower bound for the upper tail and both bounds for the lower tail.  Our main inspiration here are the works of Ganguly-Hegde \cite{GH} and Basu-Ganguly-Hammond-Hegde \cite{BGHH} which consider general last passage models. In particular, the work \cite{GH} shows how under only concavity assumptions on the limit shape, one can ``bootstrap'' a weak tail estimate of the form $\e^{ - c |s|^{\alpha}}$ to an estimate with the optimal exponents, using probabilistic and geometric techniques. Some of the techniques of \cite{GH} that we use rely   on constructions that were first completed in \cite{BGHH}.  The assumptions of \cite{GH} have been verified in only the most well-understood, integrable last passage models. Moreover, these constructions have only been carried for last passage models (the zero temperature version of polymers) and have not yet been considered for any polymer model.

The proof of the upper bound for the lower tail of \cite{GH} relies on the construction of the ``geodesic watermelon'' of \cite{BGHH}. That is, the weight of the geodesic is compared to the total weight of a large number of non-intersecting geodesics. For lower bounds, one can further restrict the paths to lie in disjoint regions of the phase space, in order to take advantage of the spatial independence of the underlying environment. Our contribution here is to adapt this construction to the polymer case, by finding estimates for non-intersecting multi-path O'Connell-Yor polymers.

An input required for this construction is a weak exponential upper bound for the lower tail. This is an assumption of \cite{GH} but does not appear in the literature for our model (the work \cite{LS} gives only estimates for the upper tail of the non-stationary model and estimates for both tails of the stationary models). The work \cite{SV} deduces variance estimates for the non-stationary model from the stationary one. By following the proof given there and inserting stronger estimates that we have derived for the stationary models using the techniques of \cite{LS}, we are able to arrive at an initial estimate of the form $\e^{ - c |s|^{3/2}}$ for an upper bound of the lower tail. 

Given this as input, we then attempt to apply the construction of \cite{BGHH} to the semi-discrete polymer case. The main super-additivity property, that $Z_{(s, m), (t, n) } \geq Z_{(s, m), (u, p)} Z_{(u, p), (t, n)}$ (here $Z_{p, q}$ is the partition function of all up-right paths from $p$ to $q$) luckily still holds and is one of the main drivers of the various proofs. However, substantial difficulties are introduced by (a) the fact that we are at positive temperature, and so $\log Z_{n, t}$ can take negative values, and (b) the semi-discrete nature of the phase space. The latter difficulty is only seen once one attempts to prove transversal fluctuation estimates and will be discussed later.

Due to the fact that $\log Z_{n, t}$ can be negative, we are forced to substantially modify the construction of \cite{BGHH}. Whereas there, some terms can be simply dropped due to the fact that a weight is always non-negative in last passage percolation, here we have to introduce a dyadic sequence of branching steps, where polymer paths split in two, and then separate from each other. This branching phase is responsible for the logarithmic loss in the range of validity  of \eqref{eqn:main-lower-upper}. 

An additional component of our work that is an input to both the upper and lower bounds for the lower tail, is handling transversal fluctuations. The work \cite{BGHH} adapts an argument of Basu, Sidoravicious and Sly \cite{BSS} which finds estimates for the geodesic weight of paths constrained to have large transversal fluctuation from the diagonal. We too adapt this argument; here the semidiscrete nature of the polymer space causes complications in the estimation of point-to-point polymer partition functions by the product of a point-to-line and line-to-point polymer partition function. However, the explicit form of the polymer partition function as well as Brownian deviation estimates allows for this sort of an estimate. The proof of the lower bound for the lower tail requires iterating this kind of estimate a number of times that grows with $n$. This is the source of the logarithmic loss in the range of validity of \eqref{eqn:main-lower-lower}. Modulo this difference, the transversal fluctuation estimates and proof of lower bound follow roughly the strategy of \cite{GH}. In particular, we rely on a version of the Harris-FKG inequality for the O'Connell-Yor polymer (we provide a proof of a version sufficient for our purposes by approximation by discrete processes in an appendix). 

One useful estimate that comes out of the treatment of transversal fluctuations that is worth separating from the rest of the paper is Corollary \ref{cor:tv} which gives,
\beq
\pp\left[ Q_{n, t} [ \TF ( \gamma ) > b n^{2/3} ] > \e^{ - c b^2 n^{1/3} }\right] \leq C \e^{ - c b^3}
\eeq
for $b \leq n^{1/3}$, where $Q_{n, t}$ denotes the polymer Gibbs measure (defined in Section \ref{sec:path}), $\gamma$ is the up-right path formed by interpreting the jumps $\{ s_i \}_i$ as the jumps of the up-right path $\gamma$ taking integer values, and $\TF ( \gamma)$ is the maximum distance of the path $\gamma$ from the straight line connecting $(0, 0)$ to $(t, n)$. In particular, this is significantly stronger than the annealed estimate that was derived for the stationary polymer in \cite{LS} using only the Rains-EJS identity and monotonicity/convexity, as the current estimate bounds the entire polymer path, instead of only the deviation at a single point, and the estimate of \cite{LS} would see the $\e^{ - c b^2 n^{1/3}}$ factor replaced by the weaker $\e^{ - c b^3}$. 

An additional wrinkle worth pointing out in the adaptation of \cite{GH} to polymer models is that the assumptions of \cite{GH} for the limit shape do not hold as written for our polymer model. For the last passage models the horizontal and vertical directions are interchangeable and so the derivative of the limit shape along the transverse direction at the diagonal vanishes; this is not the case for the O'Connell-Yor polymer. The linear correction term must therefore be accounted for when considering point-to-line or line-to-line type polymers. We carry this out by instead introducing ``compensated'' polymers; i.e., subtracting off the linear correction term.

Finally, as in \cite{GH}, the lower bound for the upper tail is a relatively straightforward consequence of super-additivity and convergence to the Tracy-Widom GUE distribution. 

It is worth mentioning that all proofs except the lower bound for the upper tail$^1$\let\thefootnote\relax\footnotetext{1. The integrable input for the lower bound for the upper tail is in fact only, roughly, that the distribution of the polymer (on the correct $n^{1/3}$ scale) is not asymptotically bounded above. Another alternative substitute would be a lower bound for the variance on the correct scale. }  make no use of integrable probability or exact formulas for the distribution of observables of the system, beyond the Burke property and stationarity (at one point we cite an estimate from \cite{SV} that is proven using the fact that Brownian LPP has the same distribution as the GUE; however the required tail estimate is a consequence of the theory of Gaussian processes and does not require this connection - see \cite{HMO}), and are probabilistic and geometric in nature. Overall, the key inputs are the Burke property of the stationary polymer, the Rains-EJS identity (which is in our setting a simple consequence of the Girsanov-Cameron-Martin formula) as well as the independence properties of the phase space combined with super-additivity of the polymer and concavity of the limit shape.

\subsection{Notational conventions}

For $a < b$ we set $[\![ a, b ] \!] := \{ m \in \zz : a \leq m \leq b \}$. For nonnegative quantities $a(i), b(i)$ depending on a parameter $i$ in an index set $\I$ (such as $n$ in the definition of the polymer) we say that $a \asymp b$ if there are $c, C>0$ so that $c a(i) \leq b(i) \leq C a(i)$ for all $i \in \I$. 

For $x = (n-m, n+m)$ we set,
\beq \label{eqn:ad-def}
\ad (x) := m ,
\eeq
(here $\ad$ stands for anti-diagonal).  Since we are dealing with up-right paths, we will often need to refer to distance between points along the diagonal and anti-diagonal axes. The diagonal distance between the points $(0, 0)$ and $(n, n) - (m, -m)$, for $|m| \leq n$ is $n$ and their anti-diagonal displacement is $|m|$. We will often say that points lying on the line $\{ (x, y) :  x + y = 2 \ell\}$ have height $\ell$. 

\subsection{Organization}

In Section \ref{sec:prelim} we collect various preliminary results from the literature about the O'Connell-Yor polymer as well as its stationary version. In Section \ref{sec:weak-lower} we establish a suboptimal upper bound on the lower tail of the form $\e^{ - c s^{3/2} }$ that is used as an a-priori input to the remainder of our paper.

The lower bound for the lower tail is carried out in Sections \ref{sec:int-int}, \ref{sec:tv} and \ref{sec:lower-bound}. In more detail, in Section \ref{sec:int-int} we establish estimates for interval-to-interval polymers. In Section \ref{sec:tv} we establish estimates for the partition function of polymers where the path is constrained to have a large transversal fluctuation. This also leads to an estimate of the quenched probability that a path has a large transversal fluctuation which is used later in the upper bound for the lower tail. In Section \ref{sec:lower-bound} we use these elements to prove the lower bound on the lower tail.

The upper bound on the lower tail takes place in Sections \ref{sec:constrained} and \ref{sec:watermelon}. In Section \ref{sec:constrained} we establish estimates on the partition function of polymers where the path is constrained to not have a large transversal fluctuation. In Section \ref{sec:watermelon} we use the watermelon construction of \cite{BGHH} to obtain the desired upper bound. 

Finally in the short Section \ref{sec:upper-tail} we obtain a lower bound for the upper tail via a short super-additivity argument.

\paragraph{Acknowledgements.} The work of B.L. is supported by an NSERC Discovery grant. B.L. thanks Amol Aggarwal and Duncan Dauvergne for helpful and illuminating discussions. The work of P.S. is partially supported by NSF grants DMS-1811093 and DMS-2154090.

\section{Preliminaries} \label{sec:prelim}

In this section we collect notation, definitions, and results from the literature useful for our work.

\subsection{Stationary and non-stationary models} \label{sec:models}

We will need to embed the O'Connell-Yor polymer in a larger family of models. First, extend $\{ B_i (s)\}_i$ to an infinite family of independent two-sided Brownian motions. We will take the convention that $B_i (0) = 0$ but in all of our definitions only increments arise and so this is irrelevant. 

  For $p, q \in \rr^2$ we use the notation $p \leq q$ to denote that the inequality holds component-wise. For $(s, m) \leq (t, n)$ we introduce the point-to-point partition function:
\beq
Z_{(s, m), (t, n) } := \int_{ s < s_m \dots < s_{n-1} < t } \e^{ \sum_{k=m}^n B_k (s_k ) - B_k (s_{k-1} ) } \d s_m \dots \d s_{n-1} ,
\eeq
where we use the convention $s_{m-1} = s$ and $s_n = t$. This convention will be used repeatedly throughout the paper without further comment when similar definitions arise. By convention we also set $Z_{(s, n), (t, n) } = \e^{ B_n (t) - B_n (s) }$. 
Then $Z_{n, t} = Z_{(0, 1), (t, n)}$. We will not use the notation $Z_{n, t}$ in the remainder of the paper. 
Note that we think of the $x$-axis as the time $t$ coordinate and the $y$-axis as the spatial integer-valued coordinate. 

We make one additional convention. If $p \leq q$ does not hold, then set $Z_{p, q} = 0$, and $\log Z_{p, q} = - \infty$.

 We will also have use for the following two-parameter version of the O'Connell-Yor polymer,
 \beq
 Z_{t, n}^{( \eta, \theta)} :=  \int_{-\infty < s_0 \dots < s_{n-1} < t } \e^{ B_0 (s_0) - \eta (s_0)_- + \theta (s_0)_+ +\sum_{k=1}^n B_k (s_k ) - B_k (s_{k-1} ) } \d s_0 \dots \d s_{n-1} .
 \eeq
 Here, $x_-=\max\{0,-x\}$ and $x_+=\max\{0,x\}$ denote the negative and positive part of $x$, respectively. We will also denote the special case $Z_{t, n}^\theta = Z_{t, n}^{( \theta, \theta)}$. In this case, $Z_{t, n}^\theta$ is stationary in a specific sense that will be used later.
 
 We define now the free energy densities by
 \beq
 f_{t,n}^\theta := t \theta - n \psi_0 ( \theta), \qquad f_{t, n} = f_{t, n}^\theta\quad \text{ with } \psi_1 ( \theta ) = t/n.
 \eeq
 Note that by elementary properties of the polygamma functions the equation $\psi_1 ( \theta ) = \kappa$ has a unique solution for any $\kappa >0$. By \cite[(2.4)]{MFSV} we have
 \beq \label{eqn:stationary-expectation}
 \ee\left[ \log Z^\theta_{t, n} \right] = f^\theta_{t, n} .
 \eeq
 The following is from \cite{LS}. 
 \bep \label{prop:old-tail}
 Let $\delta >0$ be given and let $\delta n \leq t \leq \delta^{-1} n$. There are $c, C>0$ so that
 \beq
 \pp\left[ | \log Z_{t ,n}^\theta - f_{t, n}^\theta | > s n^{1/3} \right] \leq C \e^{ - c s^{3/2} }
 \eeq
 and
 \beq
 \pp\left[  \log Z_{(0, 1), (t, n) } - f_{t, n} > s n^{1/3} \right] \leq C \e^{ - c s^{3/2}}
 \eeq
 for $0 < s < c n^{2/3}$. 
 \eep

 We also have the following from \cite[Theorem 1.1]{MFSV}
 \bel \label{lem:msfv-basic}
 Fix $\delta >0$ and let $\delta n \leq t \leq \delta^{-1} n$. There are $C,c>0$ so that
 \beq
 \ee\left[ | \log Z_{(0, 1), (t, n) } - f_{t, n} | \right] \geq cn^{1/3}
 \eeq
 and
 \beq
 \ee\left[ | \log Z_{(0, 1), (t, n) } - f_{t, n} |^2 \right] \leq C n^{2/3}
 \eeq
 \eel
 
 \bec \label{cor:basic-deviation}
 Let $\delta >0$ and assume $\delta n \leq t \leq \delta^{-1} n$. 
 There is a $c_2 >0$ and $\delta_1 > 0$ so that
 \beq
 \pp\left[ \log Z_{(0, 1), (t, n) } < f_{t, n} - c_2 n^{1/3} \right] \geq \delta_1
 \eeq
 \eec
 \proof Let $\theta$ satisfy $\psi_1 ( \theta ) = t/n$. Due to the deterministic inequality $Z_{(0, 0), (t, n) } \leq Z^{\theta}_{t, n}$  we have $\ee[ \log Z_{(0, 1), (t, n) } ] \leq \ee[ \log Z_{ t, n+1}^\theta ] = f_{t, n} + \O (1)$, where we used \eqref{eqn:stationary-expectation}. Let $X = \log Z_{(0, 1), (t, n)}$. Let $2 c_1 >0$ so that $\ee[ |X- f_{t, n} | ] \geq 2 c_1 n^{1/3}$. If $\ee[ (X - f_{t,n} )_- ] \geq c_1 n^{1/3}$ then the claim follows, since we can use the inequality
 \beq
 c_1 n^{1/3} \leq \ee[ (X - f_{t, n} )_- ] \leq c_3 n^{1/3} + \pp\left[ X - f_{t, n} < - c_3 n^{1/3}\right]^{1/2} \ee[ |X- f_{t, n} |^2 ]^{1/2}
 \eeq
 to find the desired estimate after taking, say, $c_3 = c_1/2$, after applying  Lemma \ref{lem:msfv-basic}. 
  Otherwise, assume $\ee[ (X- f_{t, n} )_+] \geq c_1 n^{1/3}$. Then, using $\ee[X] \leq f_{t, n} + C$ for some $C>0$ we have,
 \beq
 c_1 n^{1/3} \leq \ee[ (X- f_{t, n} )_+]  \leq C + \ee[ (X - \ee[X] )_+ ] = C + \ee[ (X - \ee[X])_-] \leq 2 C + \ee[ (X - f_{t, n} )_- ] ,
 \eeq
 and we conclude as before. \qed

 \subsection{Gibbs measure, polymer paths} \label{sec:path}
 
  The partition function $Z_{(s, m), (t, n)}$ is the normalization constant in the following Gibbs measure on the simplex $\{ (s_m, \dots , s_{n-1} ) \in \rr^{n-m-1} : s_m  < \dots < s_{n-1} \}$ defined by
\begin{align}
& Q_{(s, m), (t, n) } [ (s_m , \dots s_{n-1} ) \in \A ]  \notag \\
:= & \frac{1}{ Z_{(s, m), (t, n) } }   \int_{ s < s_m \dots < s_{n-1} < t } \1_{ \{(s_m , \dots s_{n-1} ) \in \A \} } \e^{ \sum_{k=m}^n B_k (s_k ) - B_k (s_{k-1} ) } \d s_m \dots \d s_{n-1}
\end{align}
for Borel $A \subseteq \rr^{n-1-m}$.  

We will interpret the times $( s_m , \dots , s_{n-1} )$ as defining the jump times of a right continuous up-right polymer path $\gamma: [s, t] \to [\![ m, n ]\!] $ uniquely defined by $\gamma (s) = k, s \in ( s_{k-1}, s_k )$. We could take $\gamma$ to be left continuous instead, but this is immaterial. 

Given some set of polymer paths $\A$ we will abuse notation and denote $Q_{(s, m), (t, n) } [ \gamma \in \A]$ as the Gibbs probability that the polymer path defined by the jump times lies in the set $\A$. We will only take very simple $\A$ so there will be no measureability concerns.

In a similar fashion we denote the Gibbs measure associated to $Z_{t, n}^\theta$ by $Q_{t, n}^\theta$ and to $Z_{t, n}^{(\eta, \theta)}$ by $Q^{(\eta, \theta)}_{t, n}$. For sets of polymer paths or jump times we will use notation,
\beq
Z_{(s, m), (t, n)} [ \gamma \in \A ] := Z_{(s, m), (t, n)} Q_{(s, m), (t, n) } [ \gamma \in \A ]
\eeq
to denote the partition function restricted to this set. Similar considerations apply to $Z_{n, t}^\theta$. Later, we will introduce several other modified or related partition functions $Z$, $\tilde{Z}$ etc., usually coming with some indices, superscripts or other decorations; they will always involve integrals over a simplex and then  notation such as $Z [ \gamma \in \A]$ or $\tilde{Z} [ \gamma \in \A] $ always means to restrict the integrals defining the partition function  at hand to the set $\A$.

\subsection{Properties of limit shape}

Note that by homogeneity, $f_{ \kappa t, \kappa n } = \kappa f_{t, n}$. We will require the following concavity of the limit shape at the point $(1, 1)$.  We introduce the following two quantities for use throughout the paper,
\beq\label{eqn: mfa-def}
\mu = f_{1, 1}, \qquad \mfa := \frac{\d}{ \d w} f_{1-w, 1+w} \vert_{w = 0}.
\eeq
\bel \label{lem:expectation-expand} Let $\eps >0$. 
There are $c, C>0$ and $\mfa \in \rr$ so that for $|w| < 1-\eps$ we have,
\beq
- C w^2 \leq f_{1-w, 1+w} - \mu - \mfa w \leq - c w^2 .
\eeq
Consequently, for $|w| <(1- \eps) n$ we have,
\beq
- C w^2 n^{-1} \leq f_{n-w, n+w}  - ( \mu n + \mfa w ) \leq - c w^2 n^{-1}.
\eeq
\eel
\proof Let $g(x, y) := f_{x,y}$ and let $\theta = \theta(x, y)$ satisfy $\psi_1 ( \theta ) = x/y$. Then, for the partial derivatives we have
\beq
g_x = \theta, \qquad g_y = - \psi_0 ( \theta) ,
\eeq
as well as,
\beq
\theta_x = \frac{1}{ \psi_2 ( \theta) y}, \qquad \theta_y = \frac{ - x}{y^2 \psi_2 ( \theta) } .
\eeq
For the Hessian of $g$ we have,
\beq
\nabla^2 g = \left( \begin{matrix} g_{xx} & g_{xy} \\ g_{yx} & g_{yy} \end{matrix} \right) = \frac{1}{y^3 \psi_2 ( \theta ) } \left( \begin{matrix} y^2 & -xy \\ -xy  & x^2 \end{matrix} \right) .
\eeq
Since $\psi_2 ( \theta) < 0$ we see that
\beq
(1, -1) \left( \nabla^2 g (x, y) \right) (1, -1)^T < - c
\eeq
for some $c>0$ and all  $\{ (x, y) \in \rr^2 : 10 > x > \eps, 10 > y > \eps \}$. The claim follows from Taylor's theorem with integral remainder. \qed

\subsection{Rains-EJS identity}
Here we state an identity derived in a more general context in \cite[Proposition 6.1]{LNS}. It is the analog for the O'Connell-Yor polymer of the identity of Rains and Emrah-Janjigian-Seppalainen for last passage percolation.
\bep \label{prop:EJS}
For any $\eta, \theta >0$ we have
\beq
\ee\left[ \exp\left( ( \eta - \theta ) \log Z_{t, n}^{ ( \eta, \theta ) } \right) \right] = \exp \left( n( \psi_{-1} ( \theta) - \psi_{-1} ( \eta ) ) - \frac{1}{2} t\left( \theta^2-\eta^2\right) \right).
\eeq
\eep

\subsection{Integer coordinates}

In many places in our work we will implicitly round quantities so that they lie on the integer lattice $\zz^2$. This usually takes place when we consider points on lines $\{ (x, y) : x +y = \ell, x, y \in \zz \}$. For example, a point $(a-b, a+b) \in \{ (x, y) : x +y = \ell, x, y \in \zz \}$ where $a$ and $b$ are not necessarily integers should be understood as the point on this line closest to $(a-b, a+b)$. This is due to the fact that these coordinates will appear in arguments of the polymer partition function, e.g., $Z_{(0, 0), (a-b, a+b)}$ which makes sense only if $a+b \in \zz$. This rounding convention does not affect proofs as the errors can be absorbed into the constants that arise in our estimates.

An additional example in which this occurs is when we divide $n$ into $k$ different segments $n/k$. For example we will want to relate $Z_{(0, 0), (n, n)}$ to $k$ copies of $Z_{(0, 0), (n/k, n/k)}$. In order to do this, one should use some combination of $Z_{(0, 0), (\lfloor n/k \rfloor, \lfloor n/k \rfloor )}$ and $Z_{(0, 0), (\lceil n/k \rceil, \lceil n/k \rceil )}$, but we will ignore this in our proofs, as the modifications are trivial and only require tedious notation.

\subsection{Rescaling}

We will prove many of our theorems only along the diagonal $Z_{(0, 0), (n, n)}$. Due to the continuus nature of the time variable, estimates for $Z_{(0, 0), (t, n)}$ and $\delta t \leq n \leq \delta^{-1} t$, for some $\delta >0$ maybe reduced to $Z_{(0, 0), (n, n)}$ by rescaling the Brownian motions by a constant order factor. Estimates throughout the work will be unchanged at the cost of adjusting constants appropriately; the limit shape $f_{t, n}$ would also of course be rescaled in some fashion but all of the properties would remain unchanged.

\section{Weak bound for lower tail} \label{sec:weak-lower}

In this section we will make use of the quantity,
\beq
e_n ( \theta, t) = t - n \psi_1 ( \theta) 
\eeq
which is the expectation of the first jump time $s_0$ with respect to the annealed measure $\ee[ Q_{n, t}^\theta[\cdot]]$, as can be seen by differentiating \eqref{eqn:stationary-expectation} wrt $\theta$.

\subsection{Jump estimates} \label{sec:jump-initial}

In this section we derive tail estimates on the first jump time $s_0$ under the annealed measure $\ee Q_{t, n}^\theta$. The work \cite{LS} derives estimates that are equivalent to an upper tail bound. A lower tail bound can be derived in much the same way. However, as the set-up in \cite{LS} at first appears slightly different from that considered here, we give all the details.

\bep \label{prop:stat-exit-1} Let $\delta >0$ and let $\delta n \leq t \leq \delta^{-1} n$. 
Let $\theta_0$ satisfy $e_n ( \theta_0 , t) = 0$. There is are $C, c, \eps >0$ so that if $\theta_0 - \eps < \eta < \theta_0$, then, 
\beq
\ee\left[ Q^{\eta}_{t, n} [ s_0 > 0   ] \right] \leq C \e^{ -c n ( \theta_0 - \eta)^3}.
\eeq
\eep
\proof We have, for any $0 < r < 1$ and $\lambda > 0$ and $\theta > \eta$,
\begin{align}
Q^{\eta}_{t, n} [ s_0 > 0 ] &\leq Q^{\eta}_{t, n} [  s_0 > 0 ]^r \leq Q^{ (\eta, \theta )}_{t, n} [ s_0 > 0 ]^r \notag \\
&= \left( \frac{ Z_{t, n}^{(\eta, \theta) } [ s_0 > 0  ] }{Z^{(\eta, \theta ) }_{t, n} }   \right)^r = \left( \frac{ Z_{t, n}^{(\lambda, \theta) } [  s_0 > 0 ] }{Z^{(\eta, \theta ) }_{t, n} }   \right)^r \notag\\
&\leq ( Z^{(\lambda, \theta)}_{t, n})^r (Z_{t, n}^{(\eta, \theta)})^{-r}
\end{align}
The first inequality follows from the fact that $x\leq x^r$ for $0 \leq x \leq 1$. The second inequality follows from the fact that
\beq
\del_y Q_{t, n}^{(x, y) } [ s_0 > 0 ] = \Cov ( \1_{ \{ s_0 > 0 \} } , (s_0)_+ ) \geq 0
\eeq
where the covariance is with respect to $Q_{t, n}^{(x, y)}$. 

Choose $4 r = \theta_0 -\eta$, and $\lambda = \theta_0$ and $\theta = \eta + 2 r$. Then we have by Proposition \ref{prop:EJS} and a Taylor expansion,
\begin{align}
 \ee\left[ ( Z^{(\lambda, \theta)}_{t, n})^r (Z_{t, n}^{(\eta, \theta)})^{-r} \right]^{2} 
 \leq~& \ee\left[ \e^{ 2 r \log Z_{t, n}^{( \theta_0, \theta ) } } \right] \times \ee\left[ \e^{ -2 r \log Z_{t, n}^{( \eta, \theta ) } } \right] \notag \\
=~& \exp \left( n ( 2 \psi_{-1} ( \theta ) - \psi_{-1} ( \theta_0) - \psi_{-1} ( \eta)  + \frac{t}{2} \left( \theta_0^2 + \eta^2 - 2 \theta^2\right) \right) \notag  \\
=~& \exp \left( 8 r^3 n \psi_2 ( \theta_0) + \O ( C r^4 n ) \right) 
\end{align}
which yields the claim. \qed

\bep \label{prop:stat-exit-2} Let $\delta >0$ and assume $\delta n \leq t \leq \delta^{-1} n$. 
Let $\theta_0$ satisfy $e_n ( \theta_0 , t) = 0$. There is are $C, c, \eps >0$ so that if $\theta_0 + \eps > \eta > \theta_0$, then, 
\beq
\ee\left[ Q^{\eta}_{t, n} [ s_0 < 0] \right] \leq C \e^{ -c n ( \eta - \theta_0)^3}.
\eeq
\eep
\proof We have, for any $\theta < \eta$ and $\lambda >0$,
\begin{align}
Q^{\eta}_{t, n} [ s_0 < 0 ]  &\leq Q^{(\theta ,\eta) }_{t, n} [  s_0 < 0 ]  = \frac{ Z^{(\theta, \eta) }_{t, n}  [  s_0 < 0]  }{ Z^{(\theta, \eta) }_{t, n} } \leq \frac{ Z_{t, n}^{ ( \theta, \lambda ) }}{ Z_{t, n}^{(\theta, \eta) }}.
\end{align}
Similar to the proof of the previous proposition, the first inequality follows from $\del_{x} Q^{(xy)}_{t, n} [ s_0 < 0 ] = - \Cov ( \1_{ \{ s_0 < 0 \} } , (s_0)_- ) \leq 0$.

Choose $4 r = \eta - \theta_0$ and $\lambda = \theta_0$ and $ \theta = \eta - 2 r$. Then, by Proposition \ref{prop:EJS} and a Taylor expansion,
\begin{align}
 & \ee \left[ \e^{ 2 r \log Z_{t, n}^{(\theta, \theta_0 ) } } \right] \ee\left[ \e^{ - 2 r \log Z_{t, n}^{(\theta, \eta)} } \right] \notag \\
= & \exp\left( n ( -2 \psi_{-1} ( \theta) + \psi_{-1} ( \eta) + \psi_{-1} ( \theta_0 ) ) + \frac{t}{2} ( - \theta_0^2 - \eta^2 + 2 \theta^2 ) \right) \notag \\
=& \exp \left( 8 r^3 n \psi_2 ( \theta_0) + \O ( C r^4 n ) \right) 
\end{align} 
This completes the proof in a similar manner to the previous result. \qed

\bec \label{cor:stat-exit}
Let $\delta >0$ and assume $\delta n \leq t \leq \delta^{-1} n$. There are $C, c>0$ so that,
\beq
\ee\left[ Q^\theta_{n, t} [ | s_0 - e_n ( \theta, t) | > s n^{2/3} ] \right] \leq C \e^{ - c s^3}
\eeq
for $0 <s \leq c n^{1/3}$. 
\eec
\proof By \cite[Remark 3.1]{SV}, we have the equality in distribution,
\beq
Q^\theta_{n, t} [ s_0 > e_n ( \theta, t) + s n^{2/3} ] \stackrel{d}{=} Q^\theta_{n, t_1} [ s_0 > 0]
\eeq
where $t_1 = t- e_n ( \theta, t) - s n^{2/3}$. We may assume that $t_1 \geq 0$ or else the claim is vacuous.  Let $\theta_0$ solve,
\beq
n \psi_1 ( \theta_0) = t_1 = t - ( t - n \psi_1 ( \theta ) ) - s n^{2/3}
\eeq
which is equivalent to,
\beq
n ( \psi_1 ( \theta_0) - \psi_1 ( \theta ) ) = -s n^{2/3}
\eeq
so that $\theta_0 - \theta \asymp s n^{-1/3}$, as long as $s \leq c n^{1/3}$, some $c>0$.  We then apply Proposition \ref{prop:stat-exit-1} to $\ee[ Q_{n, t_1}^{\theta} [ s_0 > 0 ]$ finding an estimate of $C \e^{ - c s^3}$ as long as $s \leq c n^{1/3}$, where $c>0$ is taken sufficiently small to guarantee $\theta_0 - \theta \leq \eps$ where the $\eps >0$ is from the statement of Proposition \ref{prop:stat-exit-1}.

For the other tail, we have 
\beq
Q^\theta_{n, t} [ s_0 < e_n ( \theta, t) - s n^{2/3} ] \stackrel{d}{=} Q^\theta_{n, t_2} [ s_0 < 0]
\eeq
where $t_2 = t - e_n ( \theta, t) + s n^{2/3}$. Now let $\theta_0$ solve,
\beq
n \psi_1 ( \theta_0) = t_2 = t - (t - n \psi_1 ( \theta ) ) + s n^{2/3}
\eeq
so that
\beq
n ( \psi_1 ( \theta_0 ) - \psi_1 ( \theta ) ) = s n^{2/3}
\eeq
so that $\theta - \theta_0 \asymp s n^{-1/3}$. We then apply Proposition \ref{prop:stat-exit-2} and conclude in a similar manner to the other tail. \qed

\subsection{Weak tail bound for non-stationary model}

In this section we derive a sub-optimal tail estimate for the lower tail of $Z_{(0, 1), (t, n)}$ that will serve as an input for the rest of the paper.

The proof of the following is based on the proof of \cite[Lemma 2.8]{MFSV}. Compared to that result, we have better estimates available which for various quantities that arise, allowing us to conclude a better tail than what was proven in that work.

\bep \label{prop:new-tail}
Let $\delta >0$ and assume that $\delta n \leq t \leq \delta^{-1} n$. There are $c, C>0$ so that for all $0 \leq b \leq c n^{2/3}$ we have,
\beq
\pp\left[ \log Z_{(0, 1), (t, n)} -f_{t, n} < - b n^{1/3} \right] \leq C \e^{- c b^{3/2}}.
\eeq
\eep
\proof  Let $ \theta $ satisfy $\psi_1 (\theta ) = t/n$. By Proposition \ref{prop:old-tail} it suffices to prove the estimate,
\beq
\pp \left[ \frac{ Z_{t, n}^\theta}{Z_{(0, 1), (t, n)}} \geq \e^{ b n^{1/3} } \right] \leq C \e^{ - c b^{3/2}} ,
\eeq
for some $c, C>0$. 
Compared to  \cite{MFSV}, this is an improved version of the estimate (2.40) of that paper.  In order to prove the above estimate we follow the proof of \cite[Lemma 2.8]{MFSV} inserting our better estimates where appropriate. We may assume $b \geq 1$. Let $u =  \sqrt{ b} n^{2/3}$. 
Then,
\begin{align}
 & \pp \left[ \frac{ Z_{n, t}^\theta}{Z_{(0, 1), (t, n)}} \geq \e^{ b n^{1/3} } \right] = \pp\left[ \frac{ Z^\theta_{n, t} [ |s_0 | < u ] }{ Z_{(0, 1), (t, n)} Q_{n, t}^\theta [ |s_0| < u  ] } \geq \e^{ b n^{1/3} } \right] \notag \\
\leq &\pp\left[ \frac{ Z^\theta_{n, t} [ |s_0 | < u ] }{ Z_{(0, 1), (t, n)}  } \geq \frac{1}{2} \e^{ b n^{1/3} } \right] + \pp\left[ Q_{n, t}^\theta [ |s_0| < u  ] \leq 1/2 \right]
\end{align}
For the second probability, we have 
\beq
\pp\left[ Q_{n, t}^\theta [ |s_0| < u  ] \leq 1/2 \right] = \pp\left[ Q_{n, t}^\theta [ |s_0| \geq u  ] \geq 1/2 \right] \leq C \e^{ - c u^3 n^{-2} } = C \e^{ - c b^{3/2}}
\eeq
by Corollary \ref{cor:stat-exit}. Note that we used that $e_n ( \theta, t) = 0$ by our choice of $\theta$. For the other term we estimate,
\begin{align}
 & \pp\left[ \frac{ Z^\theta_{n, t} [ |s_0 | < u ] }{ Z_{(0, 1), (t, n)} } \geq \frac{1}{2} \e^{ b n^{1/3} } \right] \notag \\
\leq & \pp\left[ \frac{ Z^\theta_{n, t} [ 0 \leq s_0  < u ] }{ Z_{(0, 1), (t, n)} } \geq \frac{1}{4} \e^{ b n^{1/3} } \right] + \pp\left[ \frac{ Z^\theta_{n, t} [ -u \leq s_0  < 0 ] }{ Z_{(0, 1), (t, n)}  } \geq \frac{1}{4} \e^{ b n^{1/3} } \right] \label{eqn:lt-a1}
\end{align}
In order to estimate the first quantity, introduce the reverse system, $B^{(r)}_0 (s) = - ( B_n (t) - B_n (t-s) )$ and $B^{(r)}_i (s) = B_{n-i} (t ) - B_{n-i} (t-s)$. Denote using the superscript $(r)$  the corresponding partition functions, Gibbs measures, etc., with respect to the reversed Brownian motions, $Z^{(r)}_{ (s, 1), (t, n)}$ and $Z_{n, t}^{\theta, (r)}$, etc. For example, $Z_{(s, 1), (t, n)} = Z^{(r)}_{(0, 0), (t-s,n-1)}$ for any $s \in (-\infty, t)$. 

Let $\nu = c_1 \sqrt{b} n^{-1/3}$ for $c_1 >0$. Let $\lambda = \theta - \nu$.  We have,
\begin{align}
\frac{ Z_{(s, 1), (t, n)}}{Z_{(0, 1), (t, n)}} &= \frac{ Z^{(r)}_{(0, 0), (t-s, n-1)} }{ Z^{(r)}_{(0, 0), (t, n-1) }} \leq \frac{ Z^{\lambda, (r)}_{t-s,n-1} [ s_0 < 0 ] }{ Z^{\lambda, (r)}_{t, n-1} [ s_0 < 0 ] } \notag \\
&= \frac{ Z^{\lambda, (r)}_{t-s,n-1} }{ Z^{\lambda, (r)}_{t,n-1} } \frac{ Q^{\lambda, (r)}_{t-s,n-1} [ s_0 < 0 ] }{ Q^{\lambda, (r)}_{t,n-1} [ s_0 < 0 ] } =: \e^{ Y^{(r)}_{n-1} (t-s, t) - \lambda s } \frac{ Q^{\lambda, (r)}_{t-s,n-1} [ s_0 < 0 ] }{ Q^{\lambda, (r)}_{t, n-1} [ s_0 < 0 ] } \notag \\
&\leq   \e^{ Y^{(r)}_{n-1} (t-s, t) - \lambda s } \frac{1}{Q^{\lambda, (r)}_{t, n-1} [ s_0 < 0 ] }
\end{align}
where in the first inequality we used the second part of the following inequality (which is \cite[(2.49)]{MFSV}),
\beq \label{eqn:mfsv-comp}
\frac{ Z^\eta_{t,n} [ s_0 > 0 ] }{ Z^\eta_{s,n} [ s_0 > 0 ] } \geq \frac{ Z_{(0, 0), (t, n)}}{Z_{(0, 0), (s, n)} } \geq \frac{ Z^\eta_{t,n} [ s_0 < 0 ] }{ Z^\eta_{s,n} [ s_0 < 0 ]  }
\eeq
which holds for any $0 < s < t $ and $\eta >0$. Here, we apply this to the reversed system, that is, adding superscripts $^{(r)}$ to the partition functions above. By definition, $Y^{(r)}_{n-1} (t-s, t) = \lambda s - \log Z^{\lambda, (r)}_{n-1, t} + \log Z^{ \lambda, (r)}_{n-1, t-s}$. 

Therefore,
\begin{align}
& \pp\left[ \frac{ Z^\theta_{t,n} [ 0 \leq s_0  < u ] }{ Z_{(0, 1), (t, n)} } \geq \frac{1}{4} \e^{ b n^{1/3} } \right] = \pp \left[ \int_0^u \e^{ - B_0 (s) + \theta s } \frac{ Z_{(s, 1), (t, n)} }{ Z_{(0, 1), (t, n)} } \d s \geq \frac{1}{4} \e^{ b n^{1/3}} \right]  \notag\\
\leq ~& \pp\left[ \int_0^t \frac{ \e^{ - B_0 (s) + ( \theta - \lambda ) s + Y^{(r)}_{n-1} (t-s, t) } }{ Q^{\lambda, (r)}_{n-1, t} [ s_0 < 0 ] } \geq \frac{1}{4} \e^{ n^{1/3} b } \right] \notag\\
\leq ~& \pp\left[ Q^{\lambda, (r)}_{n-1, t} [ s_0 < 0 ] \leq 1/2\right] + \pp\left[ \int_0^u \e^{ - B(s) + Y^{(r)}_{n-1} (t-s, t) + \nu s } \d s \geq \frac{1}{8} \e^{ n^{1/3} b } \right]
\end{align}
For the first term on the RHS, we calculate,
\begin{align}
e_{n-1} ( \lambda, t ) &\leq t - n \psi_1 ( \lambda ) + C = (t - n \psi_1 ( \theta ) ) + n( \psi_1 ( \theta) - \psi_1 ( \lambda ) ) +C \notag\\
 &\leq  - c n \nu +C \leq  - c b^{1/2} n^{2/3}
\end{align}
where the last inequality uses $\psi_2 (x) < 0$ and $b \geq 1$, and holds only for large enough $n$. We also  used $ t= n \psi_1 ( \theta)$ by definition of $\theta$. Therefore by Corollary \ref{cor:stat-exit} we have,
\begin{align}
& \pp\left[ Q^{\lambda, (r)}_{n-1, t} [ s_0 < 0 ] \leq 1/2\right] 
= \pp\left[ Q^{\lambda, (r)}_{n-1, t} [ s_0 > 0 ] \geq 1/2\right] \notag \\
\leq & \pp\left[Q^{\lambda, (r)}_{n-1, t} [ s_0 > b^{1/2} c n^{2/3} + e_{n-1} ( \lambda, t) ] \geq 1/2\right] \leq C \e^{- c b^{3/2}}.
\end{align}
Now, by the Burke property (Theorem 3.3 and Theorem 3.4 of \cite{SV}) we have that $s \mapsto Y^{(r)}_{n-1} (t-s, t)$ is a Brownian motion and by construction it is independent of $B_0 (s)$. Choose $c_1 >0$ sufficiently small so that $n^{1/3} b \geq 10 \nu u$. We therefore must bound,
\begin{align}
 & \pp \left[ \int_0^u \e^{ \sqrt{2} B_0 (s) + \nu s } \d s \geq \e^{3 \nu u } \right] 
\leq  \pp \left[ \int_{-\infty} ^u \e^{ \sqrt{2} B_0 (s) + \nu s } \d s \geq \e^{3 \nu u } \right] \notag \\
\leq & \pp\left[ \e^{ \sqrt{2} B_0 (u) + \nu u } \geq \e^{ 2 \nu u } \right] + \pp \left[ \int_{-\infty}^u \e^{ \sqrt{2} (B_0 (s) - B_0 (u) )  + \nu (s-u)} \d s \geq \e^{ \nu u } \right]
\end{align} 
The first probability is less than $C \e^{ - c \nu^2 u} \leq C \e^{ - c b^{3/2}}$. By Dufresne's identity, the integral in the second probability has the same distribution as the reciprocal of a Gamma$(\nu)$ random variable and so
\beq
\pp \left[ \mathrm{Gamma} ( \nu ) \leq \e^{- \nu u } \right] = \int_0^{ \e^{ - \nu u } } \frac{ x^{\nu-1}  \e^{ - x}}{ \Gamma (\nu) } \d x \leq \frac{ \e^{ - \nu^2 u}}{ \nu \Gamma ( \nu ) } \leq C \e^{ - \nu^2 u } \leq C \e^{ - c b^{3/2}}. 
\eeq
Collecting the above, we see that,
\beq
\pp\left[ \frac{ Z^\theta_{n, t} [ 0 \leq s_0  < u ] }{ Z_{(0, 1), (t, n)} } \geq \frac{1}{4} \e^{ b n^{1/3} } \right]  \leq C \e^{ - cb^{3/2}}
\eeq
as desired.  The second term of \eqref{eqn:lt-a1} is estimated similar to the first term.  We instead choose $\lambda = \theta + \nu$ and use the first inequality of \eqref{eqn:mfsv-comp} to obtain, 
\begin{align}
\frac{ Z^\theta_{n, t} [ - u < s_0 < 0 ] }{ Z_{(0, 1), (t, n)} } =& \int_{-u}^0 \e^{ - B_0 (s) + \theta s } \frac{ Z_{(s, 1), (t, n)} }{ Z_{(0, 1), (t, n)} } \d s \notag \\
\leq & \int_{-u}^0  \frac{ \e^{ - B_0 (s) - Y^{(r)}_{n-1} (t, t-s) - ( \theta  - \lambda ) s } }{ Q^{\lambda, (r)}_{n-1, t} [ s_0 > 0 ] } \d s 
\end{align}
Everything else is identical. \qed 

We record the above result, as well as the second estimate of Proposition \ref{prop:old-tail}, in a single corollary for easy reference.
\bec \label{cor:bad-tail}
Fix $\delta >0$ and assume $\delta n \leq t \leq \delta^{-1} n$. There are $C, c>0$ so that
\beq
\pp\left[ | \log Z_{(0, 1), (t, n) } - f_{t, n} | > s n^{1/3} \right] \leq C \e^{ -c s^{3/2}}
\eeq
for all $0 \leq s \leq c n^{2/3}$. 
\eec

\section{Interval-to-interval estimates} \label{sec:int-int}

In this section, we will consider interval-to-interval partition functions. However, the linear correction to the limit shape in Lemma \ref{lem:expectation-expand} is large and must be compensated for. We therefore do not directly consider interval-to-interval partition functions and instead consider the following modified version.

  First, for any $p \leq q \in \rr \times \zz$ we define the compensated partition function:
\beq
\tilde{Z}_{p,q} := Z_{p, q} \e^{ \mfa \ad (p-q)}.
\eeq
Here, $\mfa$ was defined in \eqref{eqn: mfa-def}, and the anti-diagonal distance operator $\ad$ was defined in \eqref{eqn:ad-def}. 

As stated above, this compensates the first order term in the correction to the free energy density of $Z_{p, q}$.  For example,
\beq
\tilde{Z}_{(m, m) + (-i, i) , (n, n) + (-j, j) } = Z_{(m, m) + (-i, i) , (n, n) + (-j, j) } \e^{ \mfa \i} \e^{ - \mfa j } .
\eeq
We will need to consider various line segment-to-line segment partition functions (or interval-to-interval). We will only consider line segments with integer coordinates parallel to the anti-diagonal. That is, we will use $\ell, \ell_i$ to denote line segments of the form,
\beq
\ell := \{ (i, i) - (j, -j ) : j \in [\![ a, b ]\!] \}
\eeq
for some $a \leq b$. Then, for two line segments $\ell_1, \ell_2$ we define,
\beq \label{eqn:def-l2l}
\tilde{Z}_{\ell_1, \ell_2} := \sum_{ p \in \ell_1, q \in \ell_2 } \tilde{Z}_{p, q} .
\eeq
By considering a point to be a line segment of one point, this definition also includes interval-to-point and point-to-interval partition functions.  Note that the quantity on the RHS may be identically $0$ if $p \leq q$ does not hold for any $(p, q) \in \ell_1 \times \ell_2$. 

In what follows, we consider two line segments, 
\beq
\ell_1 := \{ (-i, i) : |i| \leq l_1 n^{2/3} \}, \qquad \ell_2 := \{ (n, n) - (w+j, -w-j) : |j| \leq l_2 n^{2/3} \}.
\eeq
for some $w$ and $l_1, l_2 >0$.   The following is an analog of the first estimate of \cite[Proposition 3.5]{BGHH}, and is proven using a similar ``stepping back strategy.''
\bep \label{prop:int-int-upper}
Let $\ell_1$ and $\ell_2$ be as above. Assume $|w| \leq 1.5 n$. There are constants $c_1 >0$ and $s_0>0$, possibly depending on $l_1, l_2$ such that
\beq
\pp\left[ \log \tilde{Z}_{\ell_1, \ell_2} > \mu n - c_1 \frac{w^2}{n} + s n^{1/3} \right] \leq \e^{ -c_1 s^{3/2} }
\eeq
for $s_0 \leq s \leq c_1 n^{2/3}$ and $n \geq s_0$. The same estimate holds also for $\tilde{Z}_{(0, 0), \ell_2}$ and $\tilde{Z}_{\ell_1, (n-w, n+w)}$. 
\eep
\proof We do the case $l_1=l_2=1$, the general case being similar. We may assume that $|w| \leq n + 2 n^{2/3}$ or else $\tilde{Z}_{\ell_1, \ell_2 } = 0$ and the claim is trivial.  Consider the points $p = - (n, n)$ and $q = (n - w , n +w ) + (n, n)$. Let $p_*$ and $q_*$ be the points in $\ell_1$ and $\ell_2$, respectively, that satisfy
\beq
\max_{ |i| \leq n^{2/3} , |j| \leq n^{2/3}} Z_{ (-i, i), (n, n) - (w+j, -w-j) } \e^{ \mfa i } \e^{ - \mfa (w+j) }  = \tilde{Z}_{p_*, q_*} .
\eeq
Necessarily we have that $p_* \leq q_*$. 
Then,
\begin{align}
\log \tilde{Z}_{\ell_1, \ell_2} &\leq C \log(n) + \log\left(  Z_{p_*, q_*} \e^{ \mfa \ad (p_*) - \mfa \ad ( q_*) } \right)
\end{align}
and,
\begin{align}
\log\left(  Z_{p_*, q_*} \e^{ \mfa \ad (p_*) - \mfa \ad ( q_*) } \right) \leq & \log\left( Z_{p,q} \e^{ \mfa \ad (p-q) } \right) \notag \\
-& \log \left( Z_{p,p_*} \e^{ \mfa \ad (p-p_*) } \right) - \log \left( Z_{q_*, q} \e^{ \mfa \ad (q_* - q) } \right)
\end{align}
since $Z_{p, p_*} Z_{p_*, q_*} Z_{q_*, q} \leq Z_{p, q}$.  Now, $p_*$ and $q_*$ depend only on the Brownian increments $\{ B_i (s) - B_i (-i ) : -i \leq s \leq 2 n - i \}_{i }$.  The Brownian increments appearing in $Z_{p, p_*}$ can be written in terms of $\{ B_i ( -i) - B_i (s) : s \leq - i \}_i$, which are independent of the increments that $p_*$ depends on. Therefore, conditional on $p_*$, the distribution of $Z_{p, p_*}$ is simply that of a point-to-point O'Connell-Yor polymer.  A similar statement holds for $Z_{q_*, q}$. The height difference between $p$ and $p_*$ is $n$. The antidiagonal displacement between the two points is $\O ( n^{2/3})$. Therefore, by Lemma \ref{lem:expectation-expand} and Corollary \ref{cor:bad-tail} we have,
\beq
\pp\left[ | \log  Z_{p,p_*} - ( \mu n + \mfa \ad (p_* - p ) ) | >   s n^{1/3} \right] \leq \e^{ - c s^{3/2}}
\eeq
for all $n^{2/3} \geq s \geq s_0$, some $s_0 >0$, as well as a similar estimate for $\log Z_{q_*, q}$.  Therefore, if $n^{2/3} \geq s \geq s_0 +1$ and $n$ is large enough we have for any $c_1 >0$ that,
\begin{align}
\pp\left[ \log \tilde{Z}_{\ell_1, \ell_2} > \mu n - c_1 w^2 + 10 s n^{1/3} \right]  & \leq 2 \e^{ - c s^{3/2} } \notag\\
+& \pp\left[ \log \left( Z_{p, q} \e^{ \mfa \ad (p-q) }\right) > 3\mu n - c_1 w^2 + s n^{1/3} \right]
\end{align}
We have by Corollary \ref{cor:bad-tail}
\beq
\pp\left[ \log Z_{p, q} > f_{q-p} + s n^{1/3} \right] \leq \e^{- c s^{3/2}}
\eeq
for $s$ sufficiently large. On the other hand, since $| \ad (p-q) | = |w|$ and $|w| \leq n+ 2 n^{2/3}$ but the height difference of $p$ and $q$ is $3n$, we have by Lemma \ref{lem:expectation-expand} that
\beq
f_{q-p} \leq 3 \mu n + \mfa \ad (q-p) -  c_2 w^2
\eeq
some $c_2 >0$. Taking $c_1 < c_2/5$ yields the claim. \qed

We also desire a lower bound.  This is an analog of \cite[Lemma 4.4]{GH}, and is proven using a similar method. 
\bep \label{prop:int-int-lower} Let $\ell_1$ and $\ell_2$ be as above. 
Let $\delta_1 >0$ and assume $|w| \leq (1- \delta_1 ) n$. There is a $\delta_2 >0$ and $c_1 >0$ so that for all $n$ large enough,
\beq
\pp\left[ \log \tilde{Z}_{\ell_1, \ell_2} < \mu n - c_1 n^{1/3} \right] \geq \delta_2.
\eeq
\eep
We will first prove the following preliminary statement, where $l_1 = l_2 = \eps$ is taken to be small.
\bel \label{lem:int-int-lower}
Let $\delta_1 >0$ and assume $|w| \leq (1- \delta_1 ) n$ and $l_1 = l_2 = \eps >0$ is small. There is a $\delta_2 >0$, $c_1 >0$ and $\eps_0 >0$ so that for all $0 < \eps < \eps_0$ we have for $n \geq n_0 (\eps)$ that,
\beq
\pp \left[ \log \tilde{Z}_{\ell_1, \ell_2} < \mu n - c_1 n^{1/3} \right] \geq \delta_2.
\eeq
\eel
\proof Let $p  = - ( \eps^{3/2} n, \eps^{3/2} n )$ and let $q = (n-w, n+w) + ( \eps^{3/2} n , \eps^{3/2} n )$. Let $p_*$ and $q_*$ be the points maximing the summand in the definition of $\tilde{Z}_{\ell_1, \ell_2}$ so that,
\beq
\log \tilde{Z}_{\ell_1, \ell_2} \leq C \log(n) + \log \tilde{Z}_{p_*, q_*}.
\eeq
Now, the height difference of $p$ and $q$ is at least $n$, and the anti-diagonal displacement is $|w| \leq (1 - \delta_1 ) n$. The height difference of $p$ and $p_*$ is $\eps^{3/2} n$ and the anti-diagonal displacement is at most $\eps n^{2/3} \leq (1-\delta_1 ) \eps^{3/2} n$ for all $n$ large enough, depending on $\eps$. 

Similar to the proof of Proposition \ref{prop:int-int-upper} we have,
\beq
\log \tilde{Z}_{p_*, q_*} \leq \log \tilde{Z}_{p,q} - \log \tilde{Z}_{p,p_*} - \log \tilde{Z}_{q_*, q} .
\eeq
Now by Lemma \ref{lem:expectation-expand} we have,
\beq
\left| f_{p_* - p} - \mu \eps^{3/2} n - \mfa \ad (p_* - p ) \right| \leq C \eps^{1/2} n^{1/3} ,
\eeq
for some $C>0$ independent of $\eps >0$. By the independence of $p_*$ from the Brownian motion terms defining $\tilde{Z}_{p, p_*}$ we then have that for any $\delta_3 >0$, there is an $M= M(\delta_3)$ so that
\beq
\pp\left[ | \log \tilde{Z}_{p, p_*} - \mu \eps^{3/2} n | > M \eps^{1/2} n^{1/3} \right] < \delta_3
\eeq
for all $n \geq n_0 = n_0 (\eps)$. We obtain a similar estimate for $\log \tilde{Z}_{q_*, q}$.  On the other hand, by Corollary \ref{cor:basic-deviation} there is a $\delta_4 >0$ and $c_1 >0$ so that for all $n$ large enough that
\beq
\pp\left[ \log Z_{p, q} < f_{q-p} - c_1 n^{1/3} \right] > \delta_4 .
\eeq
Here we use that the fact that $|w| \leq (1- \delta_1 ) n$ implies that the anti-diagonal displacement of $p$ and $q$ is at most $(1- \delta_1 ) n$ but the height difference is at least $n$. 
By Lemma \ref{lem:expectation-expand} we have,
\beq
f_{q-p} \leq (1 + 2 \eps^{3/2} ) n \mu + \mfa \ad (q-p )
\eeq
and so
\beq
\pp\left[ \log \tilde{Z}_{p, q} < (1 +2 \eps^{3/2} ) n \mu - c_1 n^{1/3} \right] > \delta_4 
\eeq
Choose now $\delta_3 < \delta_4 /3$, which fixes $M$. Then choose $\eps >0$ so that $M \eps^{1/2} < c_1/10$. Then for all $n \geq n_0 (\eps)$ we have,
\beq
\pp\left[ \log \tilde{Z}_{p_*, q_*} < n \mu - \frac{c_1}{2} n^{1/3} \right] > \frac{\delta_4}{3}.
\eeq
This yields the claim. \qed

\noindent{\bf Proof of Proposition \ref{prop:int-int-lower}}. Choose $\eps >0$ corresponding to $\delta_1/10$ from Lemma \ref{lem:int-int-lower}.  By breaking up the intervals near $(0, 0)$ and $(n, n)$ into order $\eps^{-1}$ smaller intervals of length $\eps n^{2/3}$ we find that,
\beq
\log \tilde{Z}_{\ell_1, \ell_2} \leq C | \log \eps | + \max_{k \leq \eps^{-2}} \log \tilde{Z}_{\ell_{3,k} , \ell_{4, k} }
\eeq
for some line intervals $\ell_{3, k}$ and $\ell_{4, k}$ where the anti-diagonal displacement $w_k$ between the midpoints of $\ell_{3, k}$ and $\ell_{4, k}$ satisfies $|w_k | \leq (1- \delta_1) n + n^{2/3} \leq (1- \delta_1/2 ) n$ for $n$ large enough.  Therefore, by Lemma \ref{lem:int-int-lower} and the FKG inequality Proposition \ref{prop:fkg} we have,
\beq
\pp\left[ \max_k \tilde{Z}_{\ell_{3,k} , \ell_{4, k} } < n \mu - \frac{c_1}{2} n^{1/3} \right] \geq \prod_{k \leq \eps^{-2} } \pp\left[  \log \tilde{Z}_{\ell_{3,k} , \ell_{4, k} } < n \mu - c_1 n^{1/3} \right] \geq \delta
\eeq
for some $\delta >0$ and all $n$ large enough.  The claim now follows. \qed

Finally, we require the following point-to-long line segment estimate. 

\bep \label{prop:point-line}
Let $L = \{ (n, n) - (k, -k), |k| \leq n \}$. There is an $s_0 \geq 0$ so that for all $n^{2/3} \geq s \geq s_0$ we have,
\beq
\pp\left[ \log \tilde{Z}_{(0, 0), L} > \mu n + s n^{1/3} \right] \leq \e^{ - c s^{3/2} }
\eeq
\eep
\proof We break up $L$ into order $n^{1/3}$ line segments $\ell_i$
\beq
\ell_i := \{ (n, n) - (w_i +k, -w_i-k ) : |k| \leq n^{2/3} \}
\eeq
with $w_i := i 2 n^{2/3}$.  Then,
\beq
\log \tilde{Z}_{(0, 0), L}  \leq C \log(n) + \max_{i} \log \tilde{Z}_{(0, 0), \ell_i }.
\eeq
By Proposition \ref{prop:int-int-upper} we have for some $c_1 >0$ that for all $s \geq s_0$,
\beq
\pp\left[ \log \tilde{Z}_{(0, 0), \ell_i } > \mu n + s n^{1/3} - c_1 i^2 n^{1/3} \right] \leq \e^{ - c ( s^{3/2} + |i|^{3/2} ) } .
\eeq
The claim follows from a union bound. \qed

\section{Transverse estimates} \label{sec:tv}

The goal of the present section is to estimate the behavior of the partition function restricted to polymer paths that have a large transversal fluctuation.  This takes place over the course of several steps. In Section \ref{sec:decomp} we establish an estimate for polymer paths that have a large transversal fluctuation at their midpoint. This is accomplished by decomposing the partition function over such paths into the product of a point-to-line and line-to-point polymer, where the line has a large anti-diagonal displacement. For the purposes of the subsequence section, however, it will be necessary to establish this midpoint estimate for line-to-line polymers as well.

In Section \ref{sec:tv} we use a dyadic scheme similar to \cite{BSS} to establish the same estimate as the midpoint case as when the polymer path has a large transversal fluctuation about any point. This is Theorem \ref{thm:transversal-weight} below.  We then obtain Corollary \ref{cor:tv}, an estimate for the quenched probability that a polymer path has large transversal fluctuation.

\subsection{Decomposition} \label{sec:decomp}
Recall the notion of a polymer path $\gamma$ associated to the jump times as in Section \ref{sec:path}. 
For any $b>0$ we will let $\A_{b}^{(n)}$ denote the set of polymer paths that pass to the left of the point $(n/2, n/2) - ( bn^{2/3}, - b n^{2/3} )$. That is, $\gamma \in \A_{b}^{(n)}$ if and only if $\gamma_{n/2 - b n^{2/3} } > n/2 + b n^{2/3}$  or equivalently, $s_{n/2 + b n^{2/3} } < n/2 - b n^{2/3}$. If we are considering polymer paths defined on a rectangle $\gamma :[s, t] \to [m, n]$ that does not contain this point, we say that $\gamma \in \A_{b}^{(n)}$ iff the polymer path we get by extending $\gamma : \rr \to [ m, n]$ by setting it constant on the two intervals $(-\infty, s]$ and $[t, \infty)$ satisfies $\gamma_{n/2-bn^{2/3}} > n/2+b n^{2/3}$. 

Now for any $ a>0$, consider the line segments,
\beq
\ell^{(a)}_1 := \{ (-i, i ) : |i| \leq a n^{2/3} \}, \qquad \ell^{(a)}_2 := \{ (n, n) - (i, -i) : |i| \leq a n^{2/3} \}
\eeq
and define,
\begin{align}
\tilde{Z}^{(c)}_{n, a, b} &:= \tilde{Z}_{ \ell^{(a)}_1, \ell^{(a)}_2 } [ \A_{b+a}^{(n)} ] \notag\\
&=  \sum_{|i| \leq a n^{2/3}, |j| \leq a n^{2/3} } \tilde{Z}_{(-i, i), (n, n) - (j, -j) } [  \A_{b+a}^{(n)} ] .
\end{align} 
We now derive a decomposition of $\tilde{Z}^{(c)}_{n, a, b}$ as a product of line-to-line polymers. 
Let $p = (-i, i)$ and $q= (n, n) - (j, -j)$ be points in $\ell_1^{(a)}, \ell_2^{(a)}$, respectively.  We then decompose,
\begin{align}
&Z_{(p, q)} [ \A_{b+a}^{(n)}  ] \notag \\
=~& Z_{(p, q)} [ \A_{b+a}^{(n)}  , s_i > n -i] + Z_{(p, q)} [\A_{b+a}^{(n)}  , s_i < n-i] \notag \\
=~&  Z_{(p, q)} [ \A_{b+a}^{(n)}  , s_i > n -i] +  Z_{(p, q)} [ \A_{b+a}^{(n)} , n- (i+1) < s_i < n-i] \\
+~& Z_{(p, q)} [ \A_{b+a}^{(n)} , s_i < n-(i+1)] \notag \\
=~& Z_{(p, q)} [\A_{b+a}^{(n)}  , s_i > n -i]  +Z_{(p, q)} [\A_{b+a}^{(n)}  , n- (i+1) < s_i < n-i] \notag  \\
 +~&  Z_{(p, q) } [ \A_{b+a}^{(n)}  , s_{i+1} > n - (i+1) , s_{i} < n - (i+1) ]  
+ Z_{(p, q) } [  \A_{b+a}^{(n)}  , s_{i+1} < n - (i+1) ] \notag \\
& \dots \notag \\
=~& \sum_{k =i }^{n+j\wedge (i-1)} Z_{p, q} [\A_{b+a}^{(n)}  , s_{k-1} < n- k, s_k > n-k ] \notag \\
 +& \sum_{k =i }^{n+j\wedge i-1} Z_{p, q} [ \A_{b+a}^{(n)}  , n- (k+1) < s_k < n-k ] 
\end{align}
Let us pause to state the geometric interpretation of each of the terms appearing above. We are decomposing the partition function according to where the path crosses the line $\{ (x, y) : x + y = n \}$. The constraint $ \{ s_{k-1} < n- k, s_k > n-k \} $ implies that the path passes through the point $(n-k, k)$. The constraint $\{ n- (k+1) < s_k < n-k \}$  implies that the path crosses the line $\{ (x, y) : x +y = n \}$ in the interval $\{ (n -k, k) - (s, s) : s \in (0, 1) \}$. That is, the path jumps from level $k$ to $k+1$ in the time interval $(n-(k+1), n-k)$.

As stated above, $\gamma \in \A^{(n)}_{a+b}$ if and only if $s_{n/2+ (a+b)n^{2/3}  } < n/2 - (a+b) n^{2/3}$. 
Therefore, $Z_{p, q} [ \A_{b+a}^{(n)}  , s_{k-1} < n- k, s_k > n-k ]  = 0$ if $k \leq n/2 + (a+b) n^{2/3} $ and
\beq
Z_{p, q} [ \gamma \in \A_{b+a}^{(n)}  , s_{k-1} < n- k, s_k > n-k ] = Z_{p,(n-k, k) } Z_{(n-k, k), q}
\eeq
otherwise.

Similarly, since $s_k > n - (k+1) \implies s_{k+1} > n- (k+1)$ we see that $Z_{p, q} [ \gamma \in  \A_{b+a}^{(n)} , n- (k+1) < s_k < n-k ]  = 0$ if $k+1 \leq n/2 + (a+b) n^{2/3}$ and otherwise,
\beq
Z_{p, q} [  \A_{b+a}^{(n)} , n- (k+1) < s_k < n-k ] = \int_{n-(k+1)}^{n-k} Z_{p, (s_k , k)} Z_{(s_k, k+1), q} \d s_k .
\eeq

By summation we conclude the following.
\bep \label{prop:ZIc-decomp}
We have,
\begin{align} 
\tilde{Z}^{(c)}_{n, a, b} &= \sum_{ |i|, |j| \leq a n^{2/3} } \sum_{k> n/2  + (a+b) n^{2/3}}^{n+ j \wedge i -1 } \tilde{Z}_{(-i, i), (n-k, k) } \tilde{Z}_{(n-k, k), (n, n) - (j, -j ) } \notag\\ \\
&+ \sum_{ |i|, |j| \leq a n^{2/3} } \sum^{n+j \wedge (i-1) }_{k \geq n/2  + (a+b) n^{2/3}} \int_{n-(k+1)}^{n-k} \bigg\{  \left(  Z_{(-i, i), (s_k ,k) } \e^{\mfa (i+n/2-k)} \right) \notag\\
& \phantom{asdasdfasdfsdffsdfasf} \times  \left( Z_{(s_k, k+1), (n, n) - (j, -j) } \e^{-\mfa (j+n/2-k)}\right) \bigg\} \d s_k \label{eqn:ZIc-decomp-1}
\end{align}
\eep

We require the following lemma.
\bel \label{lem:deviation}
 There is a $c>0$ so that, for any $0 < \eps < 1$, and $r >1$ we have
\beq
\pp\left[ \sup_{ 0 \leq u \leq \eps} Z_{(s, m), (t+u, n)} > \e^{ r} Z_{(s, m), (t+\eps, n) } \right] \leq c^{-1} \e^{ -c r^2 \eps^{-1} } ,
\eeq
as well as
\beq
\pp\left[ \sup_{ 0 \leq u \leq \eps } Z_{(s-u, m), (t, n)} > \e^{ r} Z_{(s-\eps, m), (t, n)} \right] \leq c^{-1} \e^{ -c r^2 \eps^{-1} }.
\eeq
\eel
\proof We have, if $n > m$,
\begin{align} \label{eqn:deviation-main}
Z_{(s, m), (t+u, n) } &= \int_{s}^{t+u} Z_{(s, m), (w, n-1) } \e^{ B_n (t+u) - B_n (w) } \d w \notag \\
&=  \e^{ B_n (t+u) - B_n (t + \eps ) } \int_{s}^{t+u} Z_{(s, m), (w, n-1) } \e^{ B_n (t+\eps) - B_n (w) } \d w  \notag \\
& \leq  \e^{ B_n (t+u) - B_n (t + \eps ) } \int_{s}^{t+\eps} Z_{(s, m), (w, n-1) } \e^{ B_n (t+\eps) - B_n (w) } \d w \notag \\
&= \e^{ B_n (t+u) - B_n (t+\eps ) } Z_{(s, m), (t+\eps, n) }
\end{align}
The above inequality is by definition an equality in the case $m=n$. 
The first estimate then follows since,
\beq
\sup_{0 < u < \eps } B_n (t+u) - B_n (t+\eps ) \sim \sqrt{\eps} | Z|
\eeq
where $Z$ is a standard normal. The second estimate follows similarly using instead the identity
\beq
Z_{(s-u, m), (t, n) } = \int_{s-u}^t \e^{ B_m (w) - B_m (s-u) } Z_{(w, m+1), (t, n) } \d w.
\eeq
\qed

Introduce now the line segment $\ell^{(m)}$ as,
\beq
\ell^{(m)} := \{ (n/2, n/2) - (k, -k) :  2 a n \geq k \geq  (b+a) n^{2/3} \}
\eeq
Recall also the definitions of $\ell^{(a)}_1$ and $\ell^{(a)}_2$ above. 

\bep \label{prop:midpoint-2} Assume $a \geq 1$ and $1 \leq b \leq  2n$. 
There are $ C,c>0$ so that the following holds.  For any $r >1$ we have,
\beq
\log \tilde{Z}^{(c)}_{n, a, b} \leq \log \tilde{Z}_{\ell^{(a)}_1, \ell^{(m)}} + \log \tilde{Z}_{\ell^{(m)}, \ell^{(a)}_2} + C +r
\eeq
with probability at least $1 - C (a n)^3 \e^{ - c r^2}$. 
\eep
\proof We will use the identity \eqref{eqn:ZIc-decomp-1}. For the sum on the first line we have,
\begin{align}
& \sum_{k> n/2  + (a+b) n^{2/3}}^{n+ j \wedge i -1 } \tilde{Z}_{(-i, i), (n-k, k) } \tilde{Z}_{(n-k, k), (n, n) - (j, -j ) } \notag \\
& \leq \left(\sum_{k> n/2  + (a+b) n^{2/3}}^{n+ j \wedge i -1 } \tilde{Z}_{(-i, i), (n-k, k) } \right) \left( \sum_{k> n/2  + (a+b) n^{2/3}}^{n+ j \wedge i -1 } \tilde{Z}_{(n-k, k), (n, n) - (j, -j ) }  \right) \notag\\
&\leq \tilde{Z}_{\ell_1^{(a)} , \ell^{(m)} } \tilde{Z}_{ \ell^{(m)} , \ell_2^{(a)} } .
\end{align}
We now turn to the terms on the second line of \eqref{eqn:ZIc-decomp-1}. By Lemma \ref{lem:deviation} with $\eps =1$ we have
\beq
 \sup_{ n-(k+1) \leq s_k \leq n-k } Z_{(-i, i), (s_k ,k) }    Z_{(s_k, k+1), (n, n) - (j, -j) } \leq \e^{r} Z_{(-i, i),  n-k, k} Z_{( n-(k+1), k+1), (n, n) - (j, -j) }
\eeq
with probability at least $1- \e^{ - c r^2}$. The rest of the estimate follows similarly to the first line of \eqref{eqn:ZIc-decomp-1} and a union bound. \qed

\bep \label{prop:midpoint-1}
There is a $c_1 >0$ so that the following holds. 
Assume $n^{C_0} \geq a \geq 1$ for some $C_0>0$.  There is a $b_0 >0$, depending only on $C_0$ so that for all $ n \geq  b \geq b_0$ we have
\beq
\pp\left[ \log \tilde{Z}_{\ell^{(a)}_1, \ell^{(m)}} > \mu n/2 - c_1 b^2 n^{1/3} \right] \leq \e^{ - c_1 b^3} 
\eeq
and,
\beq
\pp\left[ \log \tilde{Z}_{\ell^{(m)}, \ell^{(a)}_2} > \mu n/2 - c_1 b^2 n^{1/3} \right] \leq  \e^{ - c_1 b^3} 
\eeq
\eep
\proof We prove only the first estimate, the second being similar. Let us divide $\ell^{(a)}_1$ into order $a$ line segments $\ell_{1, i}$ each of length order $n^{2/3}$ with the $i$th midpoint at the point $(-a n^{2/3}, a n^{2/3} ) + (i-1/2) (n^{2/3}, - n^{2/3} )$. Divide $\ell^{(m)}$ into at most order $n^{C_0+1}$ line segments $\ell_{2, j}$ of length order $n^{2/3}$ with midpoints $(n/2, n/2) + (a+b) (-n^{2/3}, n^{2/3} ) - (j-1/2) (n^{2/3}, n^{2/3} )$. Then,
\beq
\log \tilde{Z}_{\ell^{(a)}_1, \ell^{(m)}} \leq C \log(n) + \max_{i, j} \log \tilde{Z}_{\ell_{1,i} , \ell_{2, i} }.
\eeq
We will use a union bound to bound the max on the RHS. Now, $\tilde{Z}_{\ell_{1, i} ,\ell_{2, j} }$ is a line-to-line partition function and the anti-diagonal displacement between the midpoints of $\ell_{1, i}$ and $\ell_{2, i}$ satifies $w_{ij} \geq cn^{2/3} (i+j + b)$.  Therefore, if $b$ is sufficiently large we see by Proposition \ref{prop:int-int-upper} that (note that if $|w_{ij} | \geq 1.5 \frac{n}{2} $ then $\tilde{Z}_{\ell_{1, i} , \ell_{2, j} } = 0$)
\beq
\pp\left[ \log \tilde{Z}_{\ell_{1, i}, \ell_{2, j} } > \mu n/2 -c_2 b^2 n^{1/3} \right] \leq \e^{ - c_2 (b^3 + i^3 + j^3 ) }
\eeq
The claim now follows from a union bound. \qed

The following is the analog of \cite[Proposition 6.1]{BGHH}.

\bep \label{prop:midpoint-3}
There are $c_2 >0$ and $b_0 >0$ so that for all $b \geq b_0$ and $n^{100} \geq a \geq 1$ we have,
\beq \label{eqn:midpoint-3-est}
\pp\left[ \log \tilde{Z}^{(c)}_{n, a, b} > \mu n - c_2 b^2 n^{1/3} \right] \leq \e^{ -c_2 b^3}
\eeq
\eep
\proof We may assume that $b \leq 10 n^{1/3}$ or else $\tilde{Z}^{(c)}_{n, a, b} = 0$. Let $c_1 >0$ be the constant from Proposition \ref{prop:midpoint-1}. By Proposition \ref{prop:midpoint-2} we have with probability at least $1- C n^{1000} \e^{ - c b^4 n^{2/3} }$ that,
\beq
\log \tilde{Z}^{(c)}_{n, a, b} \leq \log \tilde{Z}_{\ell^{(a)}_1, \ell^{(m)}} + \log \tilde{Z}_{\ell^{(m)}, \ell^{(a)}_2} + c_1 b^2 n^{1/3}.
\eeq
The probability of the complementary event satisfies $n^{1000} \e^{ - c b^4 n^{2/3} } \leq \e^{ -b^3}$ for $n$ sufficiently large and $b \geq 1$. 
The claim now follows from a direct application of Proposition \ref{prop:midpoint-1}. \qed

Define now $\hat{A}^{(n)}_{b+ a}$ to be the polymer paths that pass to the right of the point $(n/2, n/2) + ( (a+b)n^{2/3}, - (a+b) n^{2/3} )$. A similar proof to that given above establishes the following.
\bep \label{prop:midpoint-4}
There are $c_2 >0$ and $b_0 >0$ so that for all $b \geq b_0$ and $n^{100} \geq a \geq 1$ we have,
\beq \label{eqn:midpoint-4-est}
\pp\left[ \log \tilde{Z}_{\ell_1^{(a)} , \ell_2^{(a)} } [ \hat{A}^{(n)}_{b+a} ] > \mu n - c_2 b^2 n^{1/3} \right] \leq \e^{ -c_2 b^3}
\eeq
\eep

\subsection{Full estimate} \label{sec:full}

For any polymer path $\gamma$ we let $\TF (\gamma)$ be its transversal fluctuation, that is, the maximal distance of the path from the diagonal. The proof of the following result follows the proof of Theorem 11.1 of \cite{BSS}.
\bet \label{thm:transversal-weight} There is a $c_1 >0$ and $b_0 , n_0>0$ so that for all $ n^{1/3} \geq b \geq b_0$ and all $n \geq n_0$ we have
\beq
\pp\left[ \log Z_{(0, 0), (n, n)} [ \TF ( \gamma ) > b n^{2/3} ] > \mu n - c_1 b^2 n^{1/3} \right] < \e^{ -c_1 b^3}.
\eeq
\eet
\proof Define the dyadic points,
\beq
S_j := \{  k 2^{-j} n : 0 \leq k \leq 2^j \}.
\eeq
Choose $j_0$ so that $2^{-j_0} n \in (0.5, 1] \times \frac{b}{10} n^{2/3}$. Define,
\beq
b_j := \frac{b}{M} \prod_{i=1}^{j-1} (1 + 2^{-i/3} ), \qquad M = 2 \cdot \prod_{i=1}^\infty (1 + 2^{-i/3} ) .
\eeq
Let $T_j$ be the set of polymer paths that intersect all of the $2^j+1$ line segments,
\beq
\ell_{k,j} := \{ ( k 2^{-j}n, k 2^{-j} n ) + (-x , x) : |x| \leq b_j n^{2/3} \} ,
\eeq
for $k=0, 1 , \dots 2^j$. 
A straightforward argument using that the paths are up-right shows that
\beq
T_{j_0} \subseteq \{ \gamma : \TF ( \gamma ) \leq b n^{2/3} \}
\eeq
and so,
\beq
\{ \gamma : \TF ( \gamma ) > b n^{2/3} \} \subseteq T_{j_0}^c = \bigcup_{j=1}^{j_0} T_j^c \cap T_{j-1}
\eeq
where $T_0$ is by definition the set of all polymer paths. Therefore,
\beq
\log Z_{(0, 0), (n, n)} [ \TF ( \gamma ) > b n^{2/3} ] \leq C \log(n) + \max_{j \leq j_0} \log Z_{(0, 0), (n, n)} [ T_j^c \cap T_{j-1} ].
\eeq
We will use a union bound to estimate the max on the RHS. First, Propositions \ref{prop:midpoint-3} and \ref{prop:midpoint-4} immediately imply that
\beq
\pp\left[ \log Z_{(0, 0), (n, n)} [ T_1^c ] > \mu n - c_1 b^2 n^{1/3} \right] \leq \e^{ - c_1 b^3}
\eeq
for some $c_1 >0$. For $1 \leq k \leq 2^{j-1}$ we let $T_{jk}^{(\pm)}$ be the polymer paths that intersect the line segments $\ell_{k-1, j-1}$ and $\ell_{k, j-1}$ and pass either above or below $\ell_{2k-1,j}$ for $\pm = +$ and $\pm = -$, respectively.  We have,
\beq
T_j^c \cap T_{j-1} \subseteq \bigcup_{k=1, \sigma \in \{ +, -\}}^{2^j-1} T_{jk}^{(\sigma)}
\eeq
Since $2^{j_0} \leq n$ we have
\beq
\log Z_{(0, 0), (n, n)} [ T_j^c \cap T_{j-1} ] \leq C \log(n) + \max_{k, \sigma \in \{ +, -\}} \log Z_{(0, 0), (n, n)} [ T_{jk}^{(\sigma)} ].
\eeq
We focus now on estimate the probability that $\log Z_{(0, 0), (n, n)}[ T_{jk}^{(+)} ]$ is large. The argument for $\log Z_{(0, 0), (n, n)}[ T_{jk}^{(-)} ]$ is similar and omitted.

We will  decompose $\log Z_{(0, 0), (n, n)} [T_{jk}^{(+)} ]$ as the product of three partition functions: (i) a point-to-interval partition function; (ii) an interval-to-interval partition function of paths constrained to have large midpoint transversal fluctuations; (iii) an interval-to-point partition function. The key point is to estimate the second partition function using Proposition \ref{prop:midpoint-3}. 

Set now $z_1 = (k-1)2^{-j+1}n$, $z_2 = k 2^{-j+1} n$ and $z_0 = (2k-1) 2^{-j} n$.  These are the coordinates of the midopints of the lines $\ell_{k-1, j-1}, \ell_{k, j-1}$ and $\ell_{2k-1, j}$, respectively. 

We require some notation to furnish our decompositions. Let $\mu_i$ be the measure on $\rr^2$ that is a sum of the $\delta$ functions at the points $\{ (z_i, z_i ) + (-m, m) : m \in \zz, |m| \leq b_{j-1} n^{2/3}  \}$ for $i=1, 2$. Let $\nu_i$ be the measure on $\rr^2$ that is a sum of $1d$ Lebesgue measures on the horizontal intervals $\{ \{ (x_i -m, x_i + m ) - (s, 0) : 0 < s < 1 \} : m \in \zz, -  b_{j-1} n^{2/3} \leq m < b_{j-1} n^{2/3}  \}$ for $i=1, 2$. Then, by a similar argument to Proposition \ref{prop:ZIc-decomp} we have (see Section \ref{sec:fd} for a complete proof)
\begin{align} \label{eqn:fd}
&Z_{(0, 0), (n, n)} [ T_{jk}^{(+)} ] = \int \int \tilde{Z}_{(0, 0), p} \tilde{Z}_{p, q} [  \A ]  \tilde{Z}_{q, (n, n)}\d \mu_1 (p) \d \mu_1 (q) \notag \\
+& \int \int\bigg\{  \left(  Z_{(0, 0), (x_1, y_1 ) } \e^{ -\mfa \ad ( (x_1, y_1 ) )} \right)\left(  Z_{(x_1, y_1 + 1), (x_2, y_2) } [ \A] \e^{ \mfa \ad ( (x_1 - x_2, y_1 - y_2 ) ) } \right) \notag \\
&\times \left( Z_{ (x_2, y_2 + 1 ), (n, n) } \e^{ \mfa \ad ( ( x_2, y_2 ) ) } \right)  \bigg\} \d \nu_1 (x_1, y_1) \d \nu_2 (x_2, y_2 ) \notag \\
+ & \int \int \tilde{Z}_{(0, 0), p} \left( Z_{p, (x, y) } [  \A ] \e^{ \mfa \ad (p - (x, y) ) }   \right) \left( Z_{(x, y+1), (n, n) } \e^{ \mfa \ad ( (x, y) ) } \right) \d \mu_1 (p) \d \nu_2 (x, y) \notag \\
+ & \int \int \left( Z_{(0, 0), (x, y) } \e^{ - \mfa \ad ( (x, y) ) } \right) \left(  Z_{(x, y+1), p} [ \A] \e^{ \mfa \ad ( (x, y) - p ) } \right) \tilde{Z}_{p, (n, n) } \d \nu_1 (x, y) \d \mu_2 (p)
\end{align}
where $\A$ is the set of polymer paths passing to the left of the point $(z_0 - b_j, z_0+ b_j )$. Let now $\ell_1$ and $\ell_2$ be the line segments,
\begin{align}
\ell_i &= \{ (z_i, z_i) + (-m, m) : |m| \leq b_{j-1} n^{2/3} \}, \qquad i=1, 2
\end{align}
Then by a similar argument to the proof of Proposition \ref{prop:midpoint-2} (the details of which appear in Appendix \ref{sec:misc-2}) we have that for any $\delta >0$
\beq \label{eqn:tv-a1}
\log Z_{(0, 0), (n, n)} [ T_{jk}^{(+)} ] \leq \delta b^2 n^{1/3} + C + \log \tilde{Z}_{(0, 0), \ell_1 } + \log \tilde{Z}_{\ell_1, \ell_2 } [ \gamma \in \A ] + \log \tilde{Z}_{\ell_2, (n, n) } ,
\eeq
with probability at least $1 - C n \e^{ - c \delta^2 b^4 n^{2/3} }$.  Now, set $r := z_2 - z_1 = 2^{-j+1} n \geq c n^{2/3}$. We see that,
\beq
\tilde{Z}_{\ell_1, \ell_2} [ \gamma \in \A ] \stackrel{d}{=} Z^{(c)}_{r, \tilde{a}, \tilde{b} }
\eeq
where $\tilde{a} = b_{j-1} n^{2/3} r^{-2/3} = 2^{2(j-1)/3} b_{j-1} \asymp 2^{2j/3} b  $ and $\tilde{b} = (b_j - b_{j-1} ) n^{2/3} r^{-2/3} \asymp b 2^{j/3}$.  Then, by Proposition \ref{prop:midpoint-3} we have (note that since $r \geq c n^{2/3}$ we have that $\tilde{a} \leq r$ for sufficiently large $n$), as long as $b$ is sufficiently large,
\beq
\pp\left[ \log \left( \tilde{Z}_{\ell_1, \ell_2} [ \gamma \in \A ]  \right) > \mu r - c_2 \tilde{b}^2 r^{1/3} \right] \leq \e^{ - c_2 \tilde{b}^3} \leq \e^{ -c b^3 2^j}
\eeq
and $c_2 \tilde{b}^2 r^{1/3} \geq c_3 b^2 2^{j/3} n^{1/3}$ for some $c_3 >0$. We take the $\delta$ introduced above in \eqref{eqn:tv-a1} to satisfy $\delta < c_3/10$.  Since $ c r \leq z_1 \leq n$ we have
\beq
\pp\left[ \log \tilde{Z}_{(0, 0), \ell_1 }  > z_1 \mu + \frac{c_3}{10} b^2 2^{j/3} n^{1/3} \right] \leq  \pp\left[ \log \tilde{Z}_{(0, 0), \ell_1 }  > z_1 \mu + \frac{c_3}{10} b^2 2^{j/3} (z_1)^{1/3} \right] \leq \e^{ - c b^3 2^{j} }
\eeq
where we applied Proposition \ref{prop:point-line} in the second inequality (note that $b^2 2^{j/3} \leq C (b/n^{1/3} ) (z_1)^{2/3}$ since $z_1 \geq r = 2^{-j+1} n$ and $2^{j_0} \leq C n^{1/3} /b$, and so the inequality is applicable). A similar estimate holds for $\log \tilde{Z}_{\ell_2, (n, n)}$. Therefore, we conclude that
\beq
\pp\left[ \log Z_{(0, 0), (n, n)}  [ T_{jk}^{(+)} ] > \mu n - \frac{c_3}{2} b^2 2^{j/3} n^{1/3} \right] \leq C n \e^{ - c b^4 n^{2/3} } + C \e^{ - c b^3 2^{j} } \leq C \e^{ -c b^3 2^{j} }
\eeq
where we used $2^j \leq n^{1/3}$ in the second inequality to simplify the estimate.  Therefore, for sufficiently large $n$ it holds for all $j$ that,
\beq
\pp\left[ \log Z_{(0, 0), (n, n) } [ T_j^c \cap T_{j-1} ] > \mu n - \frac{c_3}{3} b^2 n^{1/3} \right] \leq C 2^j \e^{ - c b^3 2^{j} } \leq C 2^{-j} \e^{ - c b^3} .
\eeq
The claim now follows. \qed

\bec \label{cor:tv} There are $c, C>0$ so that for all $b$ sufficiently large,
\beq
\pp\left[ Q_{(n, n) } [ \TF ( \gamma ) > b n^{2/3} ] > \e^{ - c b^2 n^{1/3} } \right] \leq C \e^{ -c b^3} 
\eeq
and consequently,
\beq
\ee\left[ Q_{(n, n) } [ \TF ( \gamma ) > b n^{2/3} ] \right] \leq C \e^{ -c b^3} .
\eeq
\eec
\proof With probability at least $1 - C \e^{ -c b^3}$ we have that,
\beq
Z_{(0, 0), (n, n)} [ \TF ( \gamma ) > b n^{2/3} ] \leq \e^{ \mu n -  c_1 b^2 n^{1/3} }
\eeq
and
\beq
Z_{(0, 0), (n, n) } \geq \e^{ \mu n - \frac{c_1}{2} b^2 n^{1/3} } .
\eeq
The claim follows since $Q_{(n, n)} [ \TF ( \gamma ) > b n^{2/3} ] = Z_{(0, 0), (n, n)} [ \TF ( \gamma ) > b n^{2/3} ] (Z_{(0, 0), (n, n) } )^{-1}$. \qed

\section{Lower bound for lower tail} \label{sec:lower-bound}

In this section we will prove a lower bound for the lower tail of the O'Connell-Yor polymer. For $1 \leq i , j \leq k$ we let 
\beq
v_{ij} := j \left(\frac{n}{k} , \frac{n}{k} \right) + \left( \left( i-\frac{k}{2} \right) \left( \frac{n}{k} \right)^{2/3} , -  \left( i-\frac{k}{2} \right) \left( \frac{n}{k} \right)^{2/3} \right)  ,
\eeq
and $I_{ij}$ the interval with endpoints $v_{i, j}$ and $v_{i+1, j}$. Let $L_j = \bigcup_i I_{ij}$. Let $\A$ denote the set of polymer paths that intersect every $L_j$. If a path is not in $\A$, then its transversal fluctuation is at least $c k^{1/3} n^{2/3}$. Therefore by Theorem \ref{thm:transversal-weight} we have that
\beq \label{eqn:lt-1}
\pp\left[ \log Z_{(0, 0), (n, n)} [ \A^c ] \leq \mu n - c_1 k^{2/3} n^{1/3} \right] \geq  \frac{1}{2}
\eeq
for some $c_1 >0$ and all $k$ sufficiently large.  We break up the partition function $Z_{(0, 0), (n, n)} [ \gamma \in \A] $ into a product of partition functions of $L_j$ to $L_{j+1}$ polymers. In order to do so we introduce the following measures. We let $\mu_{j, 0}$ be the measure on $\rr^2$ that is a sum of delta functions on the points of $L_i$. We let $\mu_{j, 1}$ be the measure on $\rr^2$ that is a sum of $1d$ Lebesgue measures on the horizontal intervals of the form $\{ p - (s, 0) : 0 < s <1 \}_{p \in L_j}$, except for the interval corresponding to the top-left most point of $L_j$.  Then, via  similar calculations to Proposition \ref{prop:ZIc-decomp} (see Appendix \ref{sec:sd} for a proof) we have,
\begin{align}
& Z_{(0, 0), (n, n) } [  \A ] \notag \\
= & \sum_{ \sigma \in \{0, 1\}^{k-1} } \int \bigg\{  Z_{(0, 0), (x_1, y_1 ) } \e^{ - \mfa \ad ((x_1, y_1 ) )}  \left( \prod_{j=1}^{k-2} Z_{(x_j, y_j + \sigma_j), (x_{j+1} , y_{j+1} ) } \e^{ \mfa \ad ( (x_j - x_{j+1} ) , (y_j - y_{j+1} ) ) } \right) \notag \\
& \times  Z_{(x_{k-1} , y_{k-1} + \sigma_{k-1} ), (n, n) } \e^{ \mfa \ad ( (x_{k-1} , y_{k-1} ) )} \bigg\} \prod_{j=1}^{k-1} \d \mu_{j, \sigma_j} (x_j, y_j ) . \label{eqn:lb-a2}
\end{align}
We have the estimate,
\beq \label{eqn:deviation-event}
\pp \left[ \sup_{ i, j \in \zz, |i|, |j| \leq n^{100} }  \sup_{|u| \leq n^{-9} } | B_i (j n^{-10} ) - B_i (j n^{-10} + u )  |  > 1 \right] \leq C \e^{ - c n^5}
\eeq
for some $C, c>0$.  On the complement of the event on the LHS of \eqref{eqn:deviation-event}, by the proof of Lemma \ref{lem:deviation} we see that (see Appendix \ref{sec:misc-1} for details),
\beq \label{eqn:lb-a1}
\int_{u}^{u+1} Z_{ (s, m), ( w, p ) } Z_{ (w, p+1), (t, q)} \d w \leq 10 \frac{1}{ n^{10}} \sum_{i=1}^{n^{10}} Z_{(s, m) , (u + j n^{ -10} , p ) } Z_{ (u+ (j-1) n^{-10} , p+1 ) , (t , q) }. 
\eeq
Let now $\mu_{j, 2}$ be the measure that is $n^{-10}$ times the sum of delta functions located at the points $\{ p - (m n^{-10}, 0) : p \in L_j , 0 \leq m \leq n^{10} \}$, except when $p$ is the top left point of $L_j$. That is, $\mu_{j, 2}$ is simply a discretization of $\mu_{j, 1}$ to a fine mesh. Using \eqref{eqn:lb-a1} whenever there appears $d \mu_{j, 1}$ in \eqref{eqn:lb-a2}, we have
\begin{align}
& Z_{(0, 0), (n, n) } [  \A ] \notag \\
\leq  & \sum_{ \sigma \in \{0, 1\}^{k-1} } \int \bigg\{  Z_{(0, 0), (x_1, y_1 ) } \e^{ - \mfa \ad ((x_1, y_1 ) )}  \left( \prod_{j=1}^{k-2} Z_{(x_j, y_j + \sigma_j), (x_{j+1} , y_{j+1} ) } \e^{ \mfa \ad ( (x_j - x_{j+1} ) , (y_j - y_{j+1} ) ) } \right) \notag \\
& \times  Z_{(x_{k-1} , y_{k-1} + \sigma_{k-1} ), (n, n) } \e^{ \mfa \ad ( (x_{k-1} , y_{k-1} ) )} \bigg\}\left(  \prod_{j=1}^{k-1} \d \mu_{j, 2\sigma_j} (x_j, y_j )  \right)\times C^k .
\end{align}
That is, up to an overall factor of $\O (C^k)$ we can replace the appearance of $ \d \mu_{j, 1}$ by $\d \mu_{j, 2}$. Then, using the fact that for nonnegative $f, g$ we have
\beq
\int f(x) g(x) \d \mu_{j, 2} (x) \leq C n^{10} \left( \int f(x) \d \mu_{j, 2} (x) \right) \left( \int g(x) \d \mu_{j, 2} (x) \right) ,
\eeq
we find that on the event of \eqref{eqn:deviation-event} that
\begin{align}
Z_{(0, 0), (n, n) } [ \A ] & \leq (C n)^{C k } \prod_{j=1}^{k} Z^{(j)}
\end{align}
where,
\beq
Z^{(j)} := \sum_{\sigma \in \{0, 1 \}^2} \int \tilde{Z}_{(x_1, y_1+ \sigma_1 ),( x_2, y_2 ) } \d \mu_{j-1,2 \sigma_1} (x_1, y_1) \d \mu_{j, 2 \sigma_2 } (x_2, y_2) .
\eeq
We conclude the following via the above discussion and the FKG inequality,  Proposition \ref{prop:fkg}.
\bep \label{prop:lt-1}
For any $c_2 >0$ there is a $c_3 >0$ so that if $k \leq c_3 n / ( \log(n))^3$ then,
\begin{align}
 & \pp\left[ \log Z_{(0, 0), (n, n)} \leq \mu n -  c_2 k^{2/3} n^{1/3} \right] \notag \\
\geq~& \pp\left[ \log Z_{(0, 0), (n, n) } [ \A^c] \leq \mu n - 2 c_2 k^{2/3} n^{1/3} \right] 
\times  \prod_{j=1}^k \pp\left[ \log Z^{(j)} \leq \mu (n/k) - 2 c_2 k^{-1/3} n^{1/3} \right] \notag\\
&- C \e^{ - c n^{5} }
\end{align}
\eep
We now turn to the proof of the following.
\bep \label{prop:lt-2}
There is a $c_1 >0$ so that,
\beq
\pp\left[ \log Z^{(j)} \leq \mu (n/k) - c_1 k^{-1/3} n^{1/3} \right] \geq \e^{ - C k}
\eeq
for all $k$ and $n$ large enough, satisfying $k \leq c_1 n / ( \log (n) )^3$.
\eep
\proof Let $r = n/k$. Recall that $I_{ij}$ is length $r^{2/3}$ and for each $j$ there are $k$ such intervals. We have,
\beq \label{eqn:lb-a3}
\log Z^{(j)} \leq C \log (k) + \max_{i_1,i_2} \log \hat{Z}_{I_{i_1,j-1} , I_{i_2, j} } ,
\eeq
where $\hat{Z}_{I_{i_1,j} , I_{i_2, j} }$ is the restriction of $Z^{(j)}$ to points lying near $I_{i_1, j-1}$ and $I_{i_2, j}$; that is, it involves the measures $\mu_{j-1, 0}$ and $\mu_{j, 0}$ restricted to the points in $I_{i_1, j-1}, I_{i_2, j-2}$ and the discretized intervals of the measures $\mu_{j-1, 2}$ and $\mu_{j-1, 2}$ whose right endpoint lies in $I_{i_1, j-1}$ and $I_{i_2, j-2}$ (except again, for the intervals whose right endpoint is the top left point of $I_{i_1, j-1}$ or $I_{i_2, j}$).  Note that $\hat{Z}_{I_{{i_1}, j-1} , I_{{i_2}, j} }$ is almost a line-to-line polymer, as in the definition \eqref{eqn:def-l2l}, except that we have some extra discretized horizontal segments coming from the $ \d \mu_{j, 2}$. In a moment we will replace these discretized polymers by bonafide line-to-line polymers.

From \eqref{eqn:lb-a3} and the FKG inequality, Proposition \ref{prop:fkg}, we see that for any $c_2 >0$ there is a $c_3 >0$ so that
\beq
\pp\left[ \log Z^{(j)} \leq \mu (n/k) - c_2 k^{-1/3} n^{1/3} \right] \geq \prod_{i_1=1}^k \prod_{i_2=1}^k \pp\left[ \log \hat{Z}_{I_{i_1,j-1} , I_{i_2,j}} \leq \mu (n/k) - 2 c_2 r^{1/3} \right]
\eeq
if $k \leq c_3 n / ( \log(n) )^3$.  We now wish to replace the discretized intervals by simple line-to-line polymers.

From the proof of Lemma \ref{lem:deviation} (i.e., the estimate \eqref{eqn:deviation-main} and the analog for times at the lower left endpoint) we have that for any $A \geq 1$ that,
\beq
\sup_{ \{ x_1 \in [a_1, a_1 +1] , x_2 \in [a_2-1, a_2] \}} \tilde{Z}_{(x_1, y_1 ),( x_2, y_2 ) }  \leq  \e^{A} \tilde{Z}_{ (a_1, y_1), (a_2, y_2 ) }
\eeq
with probability at least $1- \e^{ - c A^2}$. 

Therefore,
\beq \label{eqn:lt-lb-a1}
\log \hat{Z}_{ I_{i_1,j-1} , I_{i_2, j} } \leq C + A + \log \tilde{Z}_{I_{i_1,j-1} , I_{i_2, j}}
\eeq
with probability at least $1 - C r^2\e^{ - c A^2}$, the RHS defined as in \eqref{eqn:def-l2l}. 
For $|i_1-i_2| \leq \frac{1}{2} r^{1/3}$ we have from Proposition \ref{prop:int-int-lower} that,
\beq \label{eqn:lt-lb-a2}
\pp\left[ \log \tilde{Z}_{I_{i_1,j-1} , I_{i_2, j}} < \mu r - c_5 r^{1/3} \right] > \delta_2
\eeq
for some $\delta_2 >0$. Taking $A= r^{1/10}$  so that $C r^2 \e^{ - c A^2} \leq \frac{\delta_2}{2}$  for $r$ sufficiently large, we see from \eqref{eqn:lt-lb-a1} and \eqref{eqn:lt-lb-a2} that  for all $r$ large,
\beq
\pp \left[ \log \hat{Z}_{ I_{i_1,j-1} , I_{i_2, j} } < \mu r - \frac{c_5}{2} r^{1/3} \right] > \frac{\delta_2}{2} ,
\eeq
for $|i_1 - i_2 | \leq \frac{1}{2} r^{1/3}$. 
On the other hand, for $|i_1 - i_2| \geq M$, some $M>0$ we see from Proposition \ref{prop:int-int-upper} that
\beq \label{eqn:lt-lb-a3}
\pp\left[ \log \tilde{Z}_{ I_{i_1, j-1} , I_{i_2, j} } < \mu r - c_6 ( i_1 - i_2 )^2 r^{1/3} \right] \geq 1 - \e^{ - c_6 (i_1  - i_2 )^3 } .
\eeq
Taking $A = (i_1 - i_2)^2 r^{1/10}$ we see from \eqref{eqn:lt-lb-a1} and \eqref{eqn:lt-lb-a3} that for $r$ sufficiently large and after possibly increasing $M$ that
\beq 
\pp\left[ \log \hat{Z}_{I_{i_1, j-1} , I_{i_2, j} } < \mu r - \frac{c_6}{2} (i_1 - i_2 )^2 r^{1/3} \right] \geq 1 - C \e^{ - c (i_1 - i_2 )^2} .
\eeq
Therefore, we have for some $c_7 >0$ and for any $M_1 \geq M$ that
\begin{align}
& \prod_{i_1=1}^k \prod_{i_2=1}^k \pp\left[ \log \hat{Z}_{I_{i_1,j-1} , I_{i_2,j}} \leq \mu n - c_7 r^{1/3} \right] \notag \\
\geq & \left( \prod_{|i_1 - i_2 | \leq M_1 } \frac{\delta_2}{2} \right) \left( \prod_{ |i_1 - i_2| > M_1 } (1 - C \e^{ - c |i_1 - i_2 |^2 } ) \right) \geq c \left( \frac{ \delta}{2} \right)^{k M_1}
\end{align}
after setting $M_1$ possibly larger. This completes the proof. \qed

\bet \label{thm:lt-lower}
There is a $c_1 >0$ so that for any $1 \leq \theta \leq c_3 n^{2/3} ( \log (n) )^{-2}$, we have
\beq
\pp\left[ \log Z_{(0, 0), (n, n) } \leq \mu n - \theta n^{1/3} \right] \geq C \e^{ -c_1 \theta^3} .
\eeq
\eet
\proof From Propositions \ref{prop:lt-1}, \ref{prop:lt-2} and \eqref{eqn:lt-1} we see that for $k \leq c_2 n / \log(n)^3$ for some $c_2, c_1 >0$ we have,
\beq
\pp\left[ \log Z_{(0, 0), (n, n) } \leq \mu n - c_1 k^{2/3} n^{1/3} \right] \geq c \e^{ - C k^2} - C \e^{ c n^5} \geq c \e^{ - C k^2}
\eeq
where we used that $k \leq n$ in the second inequality. It remains to choose $k = C \theta^{3/2}$ for large $C>0$. \qed

\section{Constrained partition functions} \label{sec:constrained}

For any $\ell >0$ we define
\beq
Z^\ell_{(s, m), (t, n)} := Z^{\ell}_{(s, m), (t, n)} [ \TF ( \gamma) \leq \ell n^{2/3} ] . 
\eeq
The constrained partition functions will be a useful tool in proving our upper bounds on the left tail in the next section. The goal of the section is to prove the following. It is similar to \cite[Proposition 3.7]{BGHH}. 
\bep \label{prop:constrained-lower-tail}
Assume $ \delta n \leq t \leq \delta^{-1} n$. There is a $C>0, c>0$, depending on $\ell$ so that, 
\beq
\pp\left[  \log Z^\ell_{(0, 0), (t, n) } \leq f_{t, n}- u n^{1/3} \right] \leq C \e^{ - c  u}
\eeq
for all $u \geq 1$ and $n \geq C$.
\eep
\proof We break the proof into two different cases, depending on whether $u \geq n^{2/3}$ or not. First, let us assume that $u \leq C_0 n^{2/3}$. Set $J = u^{1/2}$. One can check the general inequality
\beq
Z^\ell_{(s, m), (v, p) } Z^\ell_{(v, p), (t, n)} \leq Z^\ell_{(s, m), (t, n)}.
\eeq
Therefore,
\beq
\pp\left[  \log Z^\ell_{(0, 0), (t, n)} \leq f_{t, n} - u n^{1/3} \right] \leq J \pp\left[ \log Z^\ell_{(0, 0), (t', n') } <  f_{t', n'} - J^{-2/3} u (n')^{1/3} \right]
\eeq
where we set $n' = n/J$ and $t' = t/J$. We write now,
\beq
\log Z^\ell_{(0, 0), (t', n')} = \log Z_{(0, 0), (t', n')} + \log Q_{(0, 0), (t', n')} [ \TF ( \gamma) \leq \ell J^{2/3} (n')^{2/3} ]
\eeq
We have by Corollary \ref{cor:tv} that,
(note that  $u \leq C_0 n^{2/3}$ implies $J^{2/3} \leq C (n')^{1/3}$) 
\beq
\pp\left[  \log Q_{(0, 0), (t', n')} [ \TF (\gamma) <  \ell J^{2/3} (n' )^{2/3} ]  > - C \e^{ - c u} \right] \geq 1 - C \e^{ - c u}.
\eeq

Therefore,
\begin{align}
& \pp\left[ \log Z^\ell_{(0, 0), (t' , n')} <  f_{t', n'}  - J^{-2/3} u (n')^{1/3} \right] \notag \\
\leq & \pp\left[\log Z_{(0, 0), (t' , n')} <  f_{t', n'} - J^{-2/3} u (n')^{1/3} + C \e^{- c u } \right] + C \e^{ - c u} \notag \\
\leq & C \e^{ -c u^{3/2} J^{-1} } + C \e^{ -c u} \leq C \e^{ - c u} ,
\end{align}
where we used Corollary \ref{cor:bad-tail} in the last inequality. 
This completes the proof for $u \leq C_0 n^{2/3}$. Now assume $u \geq C_0 n^{2/3}$. Take $J = C_1 n^{1/3}$. Then for $C_1 > 0$ sufficiently large, depending on $\ell>0$, we have that $Z^\ell_{(0, 0), (t J^{-1}, n J^{-1})} = Z_{(0, 0), (t J^{-1}, n J^{-1})}$. Therefore,
\beq
\pp\left[  \log Z^\ell_{(0, 0), (t, n)} \leq f_{t, n} - u n^{1/3} \right] \leq J \pp\left[ \log Z_{(0, 0), (c_1 t^{2/3}, c_2 n^{2/3} ) } < f_{c_1 t^{2/3}, c_2 n^{2/3}  } - c_3 u \right]
\eeq
By \cite[Lemma 2.9]{MFSV} we have, for $C_0$ sufficiently large,
\beq
\pp\left[ \log Z_{(0, 0), (c_1 t^{2/3}, c_2 n^{2/3} )  } < f_{c_1 t^{2/3}, c_2 n^{2/3} } - c_3 u \right] \leq C \e^{ - c' u} .
\eeq
This yields the claim. \qed

The following is an elementary consequence of standard estimates of the tail of the maximum of Brownian motion.
\bep \label{prop:sub-gaussian}
Suppose that there are $c, C>0$ so that $c \leq n, t, \ell \leq C$. Then there is a $c_1 >0$ so that,
\beq
\pp\left[ | \log Z^{\ell}_{(0, 0), (n, t)} | > u \right] \leq (c_1)^{-1} \e^{ - c_1 u^2}. 
\eeq
\eep

\section{Watermelon construction}
\label{sec:watermelon}

In this section, we show how to use the construction of \cite[Section 8]{BGHH} to get the upper bound for the lower tail. However, there are significant difficulties introduced by the fact that the log-polymer partition function can take negative values. Compared to \cite{BGHH}, we are forced to introduce a ``branching stage'' in the construction below which results in a logarithmic loss in the range of our tail bounds compared to the last passage case.

Let us take $k = 2^{N_1}$ for some $N_1 >0$ such that $k \leq c_0 n$, some $c_0 > 0$.  We begin with an informal discussion and sketch of the methodology. The basic idea is to lower bound,
\beq
\log Z_{(0, 0), (n, n)} = \frac{1}{k} \sum_{i=1}^k \log Z_{(0, 0), (n, n)} \geq \frac{1}{k} \sum_{i=1}^k \log \hat{Z}^{(i)}_{(0, 0), (n, n)}
\eeq
where $\hat{Z}^{(i)}_{(0, 0), (n, n)}$ is a carefully chosen constrained partition function. That is, $\hat{Z}^{(i)}_{(0, 0), (n, n)}$ will be an integral over polymer paths with the same Brownian increment weights as $Z_{(0, 0), (n, n)}$, however the integral will be only over paths obeying certain constraints. The constraints will be of the form that the paths have to pass through certain points in the $( t, n)$ plane and lie within a certain distance of the straight line connecting consecutive points.  There will be $k$ distinct paths/constraints, temporarily indicated by the notation $\hat{Z}^{(i)}$, and the distinct paths will spend a good amount of time in disjoint regions of phase space.

We will make repeated use of inequalities such as
\beq
Z_{(s, m), (v, p)} Z_{(v, p), (t, n)} \leq Z_{(s, m), (t, n)}
\eeq
(the LHS being interpreted as the partition function of polymer paths on $[s, t]$ starting at $m$ ending at $n$, constrained so that $\gamma_v = p$) 
and $Z^\ell_{m, n} (s, t) \leq Z_{m, n} (s, t)$.

The constraints will imply that the paths spend significant time in disjoint regions of the square $\{ (s, m ) : 0 \leq m \leq n, 0 \leq s \leq n \}$; independence will then allow for the application of concentration estimates showing that,
\beq
\pp \left[ \frac{1}{k} \sum_{i=1}^k \log \hat{Z}^{(i)}_{(0, 0), (n, n)} \leq (n - n \psi_1 ( \theta ) ) - k^{5/3} n^{1/3} \right] \leq \e^{ - c k^2}.
\eeq
We now recall our terminology that is used in order to discuss the nature of the constraints. We will be breaking up the paths into segments that are constrained to pass through points located on lines of the form $\{ x + y = 2 \ell \}$.  It is therefore convenient to use \emph{height} to refer to distance along the diagonal -- that is, points on the line $ \{ x + y = 2 \ell \}$ will be said to be at height $\ell$. A point of the form $(x + y, x-y)$ will be said to have \emph{anti-diagonal displacement} $y$. 

If the polymer paths are constrained to pass through two points $(s, m)$ and $(t, p)$ then typically we will constrain them to lie in \emph{corridors} of some \emph{width} $2 w$. That is, the polymer paths will satisfy that the maximal distance of the path $\gamma$ from the straight line connecting $(s, m)$ and $(t, n)$ will be less than $w$; that is, the paths lie within a region of width $2 w$ centered on the straight line between the points $(s, m)$ and $(t, p)$. We will say that the corridor has height $\ell$ where the point $(t-s, n-m)$ lies on the line $\{ (x, y) : x+y = 2 \ell \}$ and anti-diagonal displacement $z$ where $(t-s, n-m ) = ( \ell +z , \ell -z )$.

In general,  the corridors we consider will be of height $r$, anti-diagonal displacement $\O ( r^{2/3} )$ and width $\O ( r^{2/3} )$. That is, the anti-diagonal displacement will not be too great compared to the corridor width.

A final useful concept is the notion of separation between adjacent paths. Generically, the $k$ paths will be constrained to pass through some $k$ points $\{ p_i \}_i$ on a line $\{ x +y  =  2 \ell\}$. We will use \emph{separation} to refer to the distance along this line between consecutive points $\{p_i \}_i$.

The constraints on the polymer paths will be given as a series of five ``phases.''  We will take six heights, $\{ h_m \}_{m=0}^5$ with $h_0 = 0, h_5 = n$ and $h_1 \asymp k$, $h_2 = n/3 + \O (1)$, $h_3 = n - h_2$, $h_4 = n-h_1$.  The $m$th phase will then refer to the constraints on the polymer paths as they pass between height $h_{m-1}$ and $h_m$. The first three phases are called
\begin{enumerate}[label=(\arabic*)]
\item Branching phase
\item Separation phase
\item Middle phase
\end{enumerate}
The fourth and fifth phases are just the reverse of the separation and branching phases, respectively.

In the branching phase, all paths begin at the point $(0, 0)$ and alternately:
\begin{enumerate}[label=(\roman*)]
\item Split into two paths
\item Double the separation between consecutive paths.
\end{enumerate}
The outcome at the end of the branching phase will be $k$ paths that all have an order $1$ separation arrayed along the line $\{ x + y = 2 h_1 \}$ where $h_1 \asymp k$. Note that the branching phase necessarily must contain separation steps, to avoid all the paths clustering in a small space. Note that our procedure contains two kinds of separation, which are distinct and of a somewhat different nature: the separation that takes place during the branching phase, and the separation phase separation.  It is important not to confuse the two notions. In the branching phase, we will, for example, only seek to produce order $1$ separation.  
 
 In the separation phase, the paths will increase their separation in a dyadic fashion from $\O (1)$ to finally $\asymp n^{2/3} k^{-2/3}$ at height $h_2$. 

In the middle phase, the paths will continue along diagonal lines, maintaining the $n^{2/3} k^{-2/3}$ separation. The paths will be constrained to lie within $\O ( n^{2/3} k^{-2/3} )$ of diagonal lines so as that the weights are independent. 

\begin{figure}[htp]
    \centering
    \includegraphics[width=10cm]{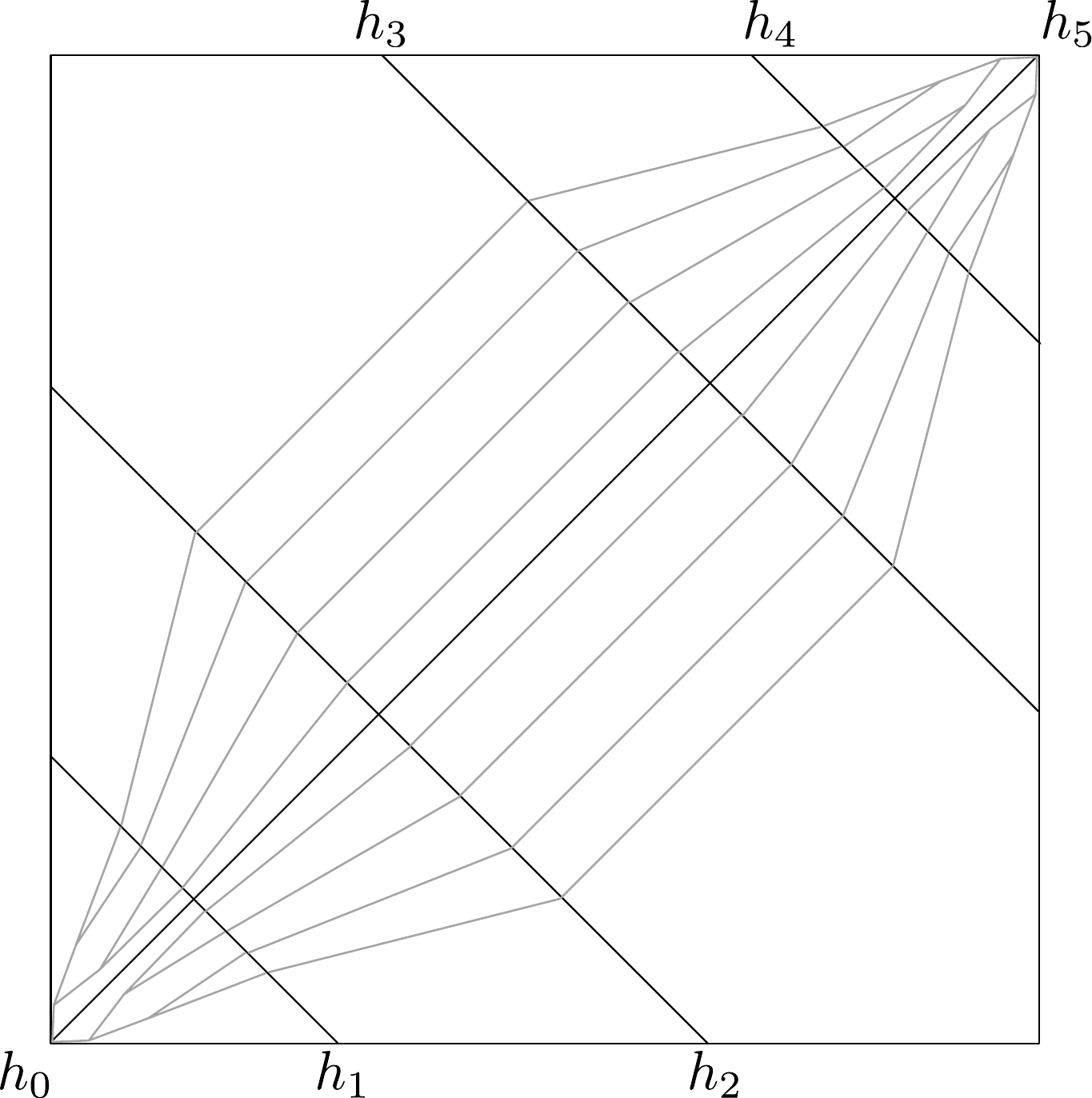}
    \caption{A schematic of the watermelon construction}
    \label{watermelon}
\end{figure}

In summary, we will estimate,
\beq
k  \log Z_{(0, 0), (n, n)} \geq \sum_{m=1}^5 \sum_{i=1}^k \log Z_{q_{i,m-1}, q_{i,m}}
\eeq
where $q_{i,m}$ are the points (to be determined) where the $i$th path intersects along the line $\{ (x , y) : x+y = 2 h_m \}$. We have $q_{i,0} = (0, 0)$ and $q_{i,5} = (n, n)$.   In the next few subsections, we will further make constraints on the paths in each of the phases, seeking lower bounds for $\sum_{i=1}^k \log Z_{q_{i,m}, q_{i,m-1} }$ for some fixed $m=1, 2, 3$ (the cases $m=4, 5$ omitted as they are similar to $j=1, 2$).

\subsection{Notational convention}

In this section we will consider points with many subscripts. With the goal of readability we will let,
\beq
Z (p, q) := Z_{p, q}
\eeq
and make similar conventions for other kinds of partition functions.  We will also denote,
\beq
f(t, n):= f_{t,n} .
\eeq

\subsection{Branching phase}

In this phase we will carry out an initial $N$ branching steps, taking us to height $h_1 \asymp k$.  That is, we will further specify constraints on the paths from $q_{i, 0}$ to $q_{i, 1}$ in order to lower bound the quantity,
\beq
\sum_{i=1}^k \log Z (q_{i, 0} , q_{i, 1} ).
\eeq
We now describe the constraints on each of the $k$ paths.

  First, every path passes from $(0, 0)$ to the vertex $(10^5, 10^5)$. Set initially, $\hat{\ell}^{(1)}_0 = 10^5$ and then $\ell^{(1)}_{j+1} = \hat{\ell}^{(1)}_j + 10^5$ for $j \geq 0$  and $\hat{\ell}^{(1)}_{j} = \ell^{(1)}_j + 10^5 \times 2^j$ for $j\geq 1$. For every $j \geq 0$,
\begin{enumerate}[label=(\roman*)]
\item Between height $\hat{\ell}^{(1)}_j$ and $\ell^{(1)}_{j+1}$ each path will branch into two paths
\item Between $\ell^{(1)}_{j+1}$ and $\hat{\ell}^{(1)}_{j+1}$, the separation between consecutive paths will increase by a factor of $2$.
\end{enumerate}
In particular, at height $\hat{\ell}^{(1)}_j$ and $\ell^{(1)}_j$ there are exactly $2^j$ distinct points that the paths intersect along the lines $\{ (x , y) : x +y = 2 \ell^{(1)}_j, 2 \hat{\ell}^{(1)}_j \}$ .  This phase ends at $h_1 = \hat{\ell}^{(1)}_{N_1}$, after there are $k = 2^{N_1}$ particles and they complete the separation step between $\ell^{(1)}_{N_1}$ and $\hat{\ell}^{(1)}_{N_1}$. 

For $j\geq 1$, at height $\ell_j^{(1)}$ we consider the $2^j$ points, $i=1, 2, \dots, 2^j$,
\begin{align}
\left( \ell_j^{(1)} - 10^4 \left( \frac{2^j+1}{2} - i \right) , \ell_j^{(1)} + 10^4 \left( \frac{2^j+1}{2} - i \right)  \right) =: \left( \ell_j^{(1)} - p_{ij}^{(1)} , \ell_j^{(1)} + p_{ij}^{(1)} \right) =: q^{(1)}_{ij}
\end{align}
and at height $\hat{\ell}^{(1)}_j$ the $2^j$ points, $i=1, 2, \dots, 2^j$,
\begin{align}
\left( \hat{\ell}_j^{(1)} -2 \times  10^4 \left( \frac{2^j+1}{2} - i \right) , \hat{\ell}_j^{(1)} + 2 \times 10^4 \left( \frac{2^j+1}{2} - i \right)  \right) =: \left( \hat{\ell}_j^{(1)} - \hat{p}_{ij}^{(1)} , \hat{\ell}_j^{(1)} + \hat{p}_{ij}^{(1)} \right) =: \hat{q}^{(1)}_{ij}.
\end{align}
At level $j$, the $k = 2^{N_1}$ paths are split into  $2^j$ equally-sized blocks of size $2^{N_1-j}$ so that the $i$th path passes through the points $q^{(1)}_{(i)_j, j}$ and $\hat{q}^{(1)}_{(i)_j, j}$ where,
\beq
(i)_j := \lceil i 2^{j-N_1} \rceil .
\eeq
By convention we set $\hat{q}^{(1)}_{i,0} = (10^5, 10^5)$. 
The above constraints are reflected in the inequality,
\begin{align}
\sum_{i=1}^k \log Z (q_{i,0}, q_{i, 1} ) & \geq  k \log Z ( (0, 0), (10^5, 10^5)  ) \label{eqn:branching-0} \\
& +   \sum_{j=1}^{N_1} \sum_{i=1}^k \log Z ( \hat{q}^{(1)}_{(i)_{j-1} , j-1 } , q^{(1)}_{(i)_j,j}  ) \label{eqn:branching-1} \\
& +  \sum_{j=1}^{N_1} \sum_{i=1}^k \log Z ( q^{(1)}_{ (i)_j ,j} , \hat{q}^{(1)}_{ (i)_j, j} ) \label{eqn:branching-2}
\end{align}
This implicitly defines the points $q_{i,1}$ as the $k$ points $\hat{q}^{(1)}_{(i)_{N_1} , N_1} = \hat{q}^{(1)}_{i, N_1}$ that have height $\hat{\ell}^{(1)}_{N_1} =h_1$. 

Note also that the diagonal distance from $q^{(1)}_{i,j}$ to $\hat{q}^{(1)}_{i,j}$ is $10^5 \times 2^j$, and the anti-diagonal displacement is at most $\pm 10^4 \times 2^j$. Therefore, the slope of the line segment connecting these two points is positive, and bounded above and away from $0$ uniformly in $j$ and $i$. The point $\hat{q}^{(1)}_{i,j-1}$ will connect to points $q^{(1)}_{2i, j}$ and $q^{(1)}_{2i-1, j}$. The diagonal distance between these points is $10^5$ and the anti-diagonal displacement is at most $\pm 10^4$, and so the lines connecting these points also has positive slope bounded above and away from $0$.

Using this decomposition, we will prove the following over the rest of this section.
\bep \label{prop:phase-1} 
There are $C, c>0$ so that,
\beq \label{eqn:phase-1-est}
\pp \left[ \sum_{i=1}^k \log Z ( q_{i,0}, q_{i, 1} )  \leq \mu k h_1 - C k^{5/3} n^{1/3} \right] \leq C \e^{ - c k^{4/3} n^{2/3} ( \log(k) )^{-1} }.
\eeq
\eep
The proof is split up into dealing with the two kinds of steps, the branching steps in \eqref{eqn:branching-1} and the separation steps in \eqref{eqn:branching-2}. 
\subsubsection{Branching steps}
We may rewrite the terms on the line \eqref{eqn:branching-1} as,
\begin{align}
\sum_{j=1}^{N_1} \sum_{i=1}^k \log Z ( \hat{q}^{(1)}_{(i)_{j-1} , j-1 } , q^{(1)}_{(i)_j,j}  ) &= \sum_{j=1}^{N_1} 2^{N_1-j} \sum_{i=1}^{2^{j-1}} \log Z ( \hat{q}^{(1)}_{i,j-1} , q^{(1)}_{2i-1,j} ) +  \log Z ( \hat{q}^{(1)}_{i,j-1} , q^{(1)}_{2i,j} ) \notag \\
&\geq \sum_{j=1}^{N_1} 2^{N_1-j} \sum_{i=1}^{2^{j-1}} \log \hat{Z}^{(1)} ( \hat{q}^{(1)}_{i,j-1} , q^{(1)}_{2i-1,j} ) +  \log \hat{Z}^{(1)} ( \hat{q}^{(1)}_{i,j-1} , q^{(1)}_{2i,j} ) \notag \\
&=: \sum_{j=1}^{N_1} 2^{N_1 - j} \sum_{i=1}^{2^{j-1}} Y^{(1)}_{ij}
\end{align}
where we define $\hat{Z}^{(1)} ( x, y)$ as the partition function of polymer paths from $x$ to $y$ constrained to lie within $10^3$ of the straight line connecting $x$ to $y$, and
\beq
Y_{ij} :=  \log \hat{Z}^{(1)} ( \hat{q}^{(1)}_{i,j-1} , q^{(1)}_{2i-1,j} ) +  \log \hat{Z}^{(1)} ( \hat{q}^{(1)}_{i,j-1} , q^{(1)}_{2i,j} )
\eeq
 Now, the collection $\{Y_{ij} \}_{ij}$ are mutually independent random variables and by Proposition \ref{prop:sub-gaussian} we have (since the width and height of the corridors involved are of constant order),
 \beq \label{eqn:sub-g-1}
 \pp\left[ | Y_{ij} | > u \right] \leq C \e^{ - c u^2}
 \eeq
for some $C, c >0$ and all $u \geq u_0$. 
The following follows from standard sub-Gaussian concentration results (see, e.g., \cite[Section 2.5]{V}).
\bel We have that,
\beq \label{eqn:branching-1-est}
\pp\left[\left| 2^{N_1} \sum_{j=1}^{N_1} 2^{-j} \sum_{i=1}^{2^{j-1}} Y_{ij} \right| > k^{5/3} n^{1/3} \right] \leq 2 \e^{ - c k^{4/3} n^{2/3} }.
\eeq
\eel
\proof This estimate follows from a direct application of Proposition \ref{prop:sub-gauss}, with $G_{ij} = Y_{ij}$ and $a_{ij} = 2^{-j}$. The estimate \eqref{eqn:sub-g-1} guarantees that the hypotheses are fulfilled. We calculate,
\beq
\| a \|_1 = \sum_{j=1}^{N_1} \sum_{i=1}^{2^{j-1}} 2^{-j} = \frac{N_1}{2} \leq C \log (k)
\eeq
and
\beq
\| a \|_2^2 = \sum_{j=1}^{N_1} \sum_{i=1}^{2^{j-1}} ( 2^{-j} )^2 \leq \sum_{j=1}^{N_1} 2^{-j} \leq C .
\eeq
Therefore, Proposition \ref{prop:sub-gauss} gives the estimate,
\beq
\pp\left[ \left| \sum_{j=1}^{N_1} \sum_{i=1}^{2^{j-1} } 2^{-j} Y_{ij} \right| \geq  C ( \log (k) + t ) \right] \leq C \e^{ - c t^2} .
\eeq
The claim now follows by taking $t = c_1 k^{2/3} n^{1/3}$ for some small $c_1 >0$. \qed

\subsubsection{Order 1 separation steps}

We begin by rewriting the terms \eqref{eqn:branching-2} as,
\begin{align}
\sum_{j=1}^{N_1} \sum_{i=1}^k \log Z ( q^{(1)}_{ (i)_j ,j} , \hat{q}^{(1)}_{ (i)_j, j} ) = 2^{N_1} \sum_{j=1}^{N_1} \frac{1}{2^j} \sum_{i=1}^{2^j} \log Z ( q^{(1)}_{ i j} , \hat{q}^{(1)}_{ i j} ) 
\end{align}
The height difference between the points $q^{(1)}_{ij}$ and $\hat{q}^{(1)}_{ij}$ is order $2^j$. The anti-diagonal displacement between the two points is as much as order $2^j$ for $i$ close to $1$ or $2^j$. In order to obtain corridors with bounded aspect ratios, we therefore split the paths passing between the levels $\ell^{(1)}_j$ and $\hat{\ell}^{(1)}_j$ into order $2^j$ further sub-levels.

For $s=0, 1, \dots 2^j$ define $\hat{\ell}^{(1)}_{j,s} = \ell_j^{(1)} + 10^5 \times s$.  Then, let $\hat{q}^{(1)}_{ij,s}$ be the point that is the intersection of the line $\{ (x , y) :  x + y = 2 \hat{\ell}^{(1)}_{j,s} \}$ and the straight line segment connection $q^{(1)}_{ij}$ and $\hat{q}^{(1)}_{ij}$. Then,
\begin{align}
2^{N_1} \sum_{j=1}^{N_1} \frac{1}{2^j} \sum_{i=1}^{2^j} \log Z ( q^{(1)}_{ i j} , \hat{q}^{(1)}_{ i j} ) &\geq 2^{N_1} \sum_{j=1}^{N_1} \frac{1}{2^j} \sum_{i=1}^{2^j} \sum_{s=1}^{2^j}  \log Z ( \hat{q}^{(1)}_{ij,s-1} ,\hat{q}^{(1)}_{ij,s} ) \notag\\
&\geq  2^{N_1} \sum_{j=1}^{N_1} \frac{1}{2^j} \sum_{i=1}^{2^j} \sum_{s=1}^{2^j}  \log \hat{Z}^{(1)} ( \hat{q}^{(1)}_{ij,s-1} ,\hat{q}^{(1)}_{ij,s} ) 
\end{align}
where again, $\hat{Z}^{(1)} (x, y)$ denotes the polymer partition function of paths from $x$ to $y$ staying within $10^3$ of the straight line connecting these points. Due to these constraints and the separation between consecutive points at each height $\hat{\ell}^{(1)}_{j,s}$ we see that the collection $\{  \log \hat{Z}^{(1)} ( \hat{q}^{(1)}_{ij,s-1} ,\hat{q}^{(1)}_{ij,s} )  \}_{ijs}$ are a family of independent random variables obeying,
\beq \label{eqn:sub-gauss-2}
\pp\left[ |  \log \hat{Z}^{(1)} ( \hat{q}^{(1)}_{ij,s-1} ,\hat{q}^{(1)}_{ij,s} ) | > u \right] \leq 2 \e^{ - c u^2} ,
\eeq
for some $c>0$, where we applied Proposition \ref{prop:sub-gaussian}. Therefore, again by sub-Gaussian concentration (see \cite[Section 2.5]{V}), we derive the following.
\bel There is a $C_1 > 0 $ so that,
\beq \label{eqn:branching-2-est}
\pp\left[ \left| 2^{N_1} \sum_{j=1}^{N_1} \frac{1}{2^j} \sum_{i=1}^{2^j} \sum_{s=1}^{2^j}  \log \hat{Z}^{(1)} ( \hat{q}^{(1)}_{ij,s-1} ,\hat{q}^{(1)}_{ij,s} ) \right| > C_1 k^{5/3} n^{1/3} \right] \leq  2 \e^{ - c k^{4/3} n^{2/3} (\log (k ) )^{-1} } ,
\eeq
\eel
\proof  This estimate follows by a direct application of Proposition \ref{prop:sub-gauss}. Let $G_{ijs} =  \log \hat{Z}^{(1)} ( \hat{q}^{(1)}_{ij,s-1} ,\hat{q}^{(1)}_{ij,s} )$ and $a_{ijs} = 2^{-j}$. The estimate \eqref{eqn:sub-gauss-2} guarantees that the hypotheses are fulfilled. We calculate,
\beq
\| a \|_1 = \sum_{j=1}^{N_1} \sum_{i=1}^{2^j} \sum_{s=1}^{2^j} 2^{-j} = \sum_{j=1}^{N_1} 2^j \leq C 2^{N_1} \leq C k
\eeq
where we used that $2^{N_1} = k$. We also calculate,
\beq
\| a \|_2^2 =\sum_{j=1}^{N_1} \sum_{i=1}^{2^j} \sum_{s=1}^{2^j} ( 2^{-j} )^2 \leq  \sum_{j=1}^{N_1} C = C N_1 \leq C \log (k).
\eeq
Therefore, from Proposition \ref{prop:sub-gauss} we conclude the estimate,
\beq
\pp\left[ \left| \sum_{j=1}^{N_1} \frac{1}{2^j} \sum_{i=1}^{2^j} \sum_{s=1}^{2^j}  \log \hat{Z}^{(1)} ( \hat{q}^{(1)}_{ij,s-1} ,\hat{q}^{(1)}_{ij,s} ) \right| > C (k + t ) \right] \leq C \exp \left( - \frac{ c t^2}{ \log (k) } \right)
\eeq
for all $t >0$. The claim follows from the choice of $t = k^{2/3} n^{1/3}$ and the fact that $k \leq n$. \qed

\subsubsection{Proof of Proposition \ref{prop:phase-1}}

First, note that since $h_1 \leq C k \leq C n$, the term $\mu k h_1$ appearing in \eqref{eqn:phase-1-est} can be absorbed into the term $k^{5/3} n^{1/3}$ at the expense of changing the constants.  Now, in order to complete the proof of the Proposition, we use the lower bound of $\sum_i \log Z (q_{i,0} , q_{i, 1} )$ in terms of the three terms \eqref{eqn:branching-0}, \eqref{eqn:branching-1} and \eqref{eqn:branching-2}. First, for the term \eqref{eqn:branching-0}, we have by Proposition \ref{prop:sub-gaussian} that,
\beq
\pp \left[ k \log Z ((0, 0) , (10^5, 10^5) ) < - C k^{5/3} n^{1/3} \right] \leq C \e^{ - k^{4/3} n^{2/3} } .
\eeq
Next, the terms \eqref{eqn:branching-1} and \eqref{eqn:branching-2} are handled by the estimates \eqref{eqn:branching-1-est} and \eqref{eqn:branching-2-est} that were proved over the course of the previous two subsections. \qed

\subsection{Separation phase} 
In this phase, set first $\ell^{(2)}_0 = h_1$, the endpoint of the previous phase. We then inductively define,
\beq
\ell^{(2)}_j = \ell^{(2)}_{j-1} + 10^5 (2^j)^{3/2} k.
\eeq
There will be $N_2$ levels where $2^{N_2} = 10^{-10} n^{2/3} k^{-2/3}(1+ \O (1))$.  Note that with this choice and the assumption $k \leq c_1 n$ some small $c_1 >0$ we have that,
\beq
c n \leq h_2 \leq n/3
\eeq
 The $i$th curve will intersect level $\ell^{(2)}_j$ at the point,
\beq
\left( \ell^{(2)}_j - p^{(2)}_{ij} , \ell^{(2)}_j + p^{(2)}_{ij}  \right) =: q^{(2)}_{ij}
\eeq
where,
\beq
p^{(2)}_{ij} = p^{(2)}_{i,j-1} + 2^j \left( \frac{k+1}{2} - i \right) ,
\eeq
and $p^{(2)}_{i0} = 2 \times 10^4 \left( \frac{k+1}{2} - i \right)$ (i.e., the last points from the previous phase). Note that,
\beq
p^{(2)}_{ij} = - p^{(2)}_{(k+1-i),j} .
\eeq
Note also that the separation between consecutive points on level $j$ is,
\beq
C 2^j \geq \mathrm{sep}_j := p^{(2)}_{i,j} - p^{(2)}_{i-1,j} \geq 2^j + 10^5.
\eeq
The antidiagonal displacement between the points that the $i$th path crosses on consecutive levels (i.e., between points $q^{(2)}_{i,j-1}$ and $q^{(2)}_{ij}$) is as much as,
\beq
p^{(2)}_{1j} - p^{(2)}_{1,j-1} = 2^j k ,
\eeq
which is larger than the separation by a factor of $k$ (note also that the points $q^{(2)}_{i,j-1}$ and $q^{(2)}_{i,j}$ have diagonal separation $10^5 (2^j)^{3/2} k$ and so the slope of the straight line connecting this points is positive and bounded above and away from $0$). We therefore split the path between consecutive levels $\ell^{(2)}_{j-1}$ and $\ell^{(2)}_j$ into $k$ further sublevels,
\beq
\ell^{(2)}_{j,s} = \ell^{(2)}_{j-1} + 10^5 (2^j)^{3/2} s
\eeq
for $s=0, \dots, k$. We then let 
\beq
q^{(2)}_{ijs} = \left( \ell^{(2)}_{js} - p^{(2)}_{ijs} , \ell^{(2)}_{js} + p^{(2)}_{ijs} \right)
\eeq
be the point that is on the intersection of the straight line connecting $q^{(2)}_{ij}$ and $q_{i,j+1}^{(2)}$ and the line $\{ (x, y) : x+ y = 2 \ell^{(2)}_{j,s} \}$. Note that,
\beq
p^{(2)}_{ijs} = - p^{(2)}_{k+1-i,j,s} + \O (1)
\eeq
the $\O(1)$ coming from the effect of the integer lattice.  We now write,
\begin{align}
\sum_{i=1}^k \log Z_{q_{i,1} , q_{i,2} } &\geq \sum_{j=1}^{N_2} \sum_{i=1}^k \sum_{s=1}^k \log Z ( q^{(2)}_{ij,s-1} , q^{(2)}_{ij,s} ) \notag \\
&\geq \sum_{j=1}^{N_2} \sum_{i=1}^k \sum_{s=1}^k \log \hat{Z}^{(2),j} ( q^{(2)}_{ij,s-1} , q^{(2)}_{ij,s} ) 
\end{align}
where we define $\hat{Z}^{(2),j}  (x, y)$ to be the partition function of polymer paths from $x$ to $y$ restricted to the corrider of width $2^j + 10^4$ around the straight line connecting $x$ to $y$. 
\bel For some $c, C>0$ we have,
\beq
\pp\left[ \sum_{j=1}^{N_2} \sum_{i=1}^k \sum_{s=1}^k \left(  \log \hat{Z}^{(2),j} ( q^{(2)}_{ij,s-1} , q^{(2)}_{ij,s} )  - f ( q_{ij,s}^{(2)} - q_{ij,s-1}^{(2)} )\right) \leq - C k^{5/3} n^{1/3} \right] \leq \e^{ - c k^2}
\eeq
\eel
\proof Let $Y_{ijs} := f ( q_{ij,s}^{(2)} - q_{ij,s-1}^{(2)} ) -\log \hat{Z}^{(2),j} ( q^{(2)}_{ij,s-1} , q^{(2)}_{ij,s} ) $.  The height and width of the rectangle with opposite vertices $q_{ij,s}$ and $q_{ij,s-1}$ and sides parallel to the coordinate axes are both order $(2^j)^{3/2}$ by our earlier discussion of the diagonal and anti-diagonal separation of $q_{i,j-1}$ and $q_{i, j}$.  
With $r_j = (2^j)^{3/2}$ we see that the polymer paths are restricted to lie in a corridor of width $r_j^{2/3}$ and so Proposition \ref{prop:constrained-lower-tail} is applicable if $j$ is sufficiently large. Therefore,
\beq
\pp\left[ Y_{ijs} \geq u C (r_j)^{1/3} \right] \leq \e^{ - u}
\eeq
for $u \geq u_0$ and $j$ sufficiently large. For smaller $j$ we instead can apply Proposition \ref{prop:sub-gaussian} to arrive at the same estimate.

We now wish to apply Proposition \ref{prop:conc}. Note that the family of random variables $Y_{ijs}$ are independent because the separation between points on level $j$ is at least $2^j + 10^5$, as discussed above, and the paths are restricted to lie in corridors of width $2^j+10^4$. 

We therefore may apply Proposition \ref{prop:conc} with $a_{ijs}^{-1} := C (r_j)^{1/3}$. Then,
\beq
\nu := \sum_{ijs} \frac{1}{a_{ijs}} = C k^2 \sum_{j=1}^{N_2} (r_j)^{1/3} = C k^2  \sum_{j=1}^{N_2} (2^j)^{1/2} \asymp k^{5/3} n^{1/3}
\eeq
by the choice $2^{N_2} \asymp n^{2/3} k^{-2/3}$, and
\beq
\min_{ijs} a_{ijs} \geq c k^{1/3} n^{-1/3}. 
\eeq
Therefore, Proposition \ref{prop:conc} implies,
\beq
\pp\left[ \sum_{ijs} Y_{ijs} \geq C k^{5/3} n^{1/3} \right] \leq \e^{ - c k^2}
\eeq
as desired. \qed

\bel
We have,
\beq
\left| (h_2 - h_1) \mu k - \sum_{j=1}^{N_2} \sum_{i=1}^k \sum_{s=1}^k f ( q_{ij,s}^{(2)} - q_{ij,s-1}^{(2)} )  \right| \leq C k^{5/3} n^{1/3} .
\eeq
\eel
\proof Writing $h_2 - h_1 = \sum_{j,s} ( \ell^{(2)}_{j,s } - \ell^{(2)}_{j,s-1})$ we have,
\begin{align}
& \sum_{j=1}^{N_2} \sum_{i=1}^k \sum_{s=1}^k f ( q_{ij,s}^{(2)} - q_{ij,s-1}^{(2)} )  - (h_2 - h_1) \mu k \notag \\
= & \sum_{j=1}^{N_2} \sum_{s,i} \left( f ( q_{ij,s}^{(2)} - q_{ij,s-1}^{(2)} ) - \mu ( \ell^{(2)}_{i,s } - \ell^{(2)}_{j,s-1}) \right) 
\end{align}
Let $r_{j,s} := ( \ell^{(2)}_{i,s } - \ell^{(2)}_{j,s-1}) $ and
\beq
\hat{p}_{i,j,s} = p^{(2)}_{ijs} - p^{(2)}_{ij,s-1}
\eeq
so that,
\beq
q_{ij,s}^{(2)} - q_{ij,s-1}^{(2)} = \left( r_{j,s} - \hat{p}_{ijs} , r_{j,s} + \hat{p}_{ijs} \right).
\eeq
Now we have 
\beq \label{eqn:phase-2-a1}
\hat{p}_{ijs} = - \hat{p}_{k+1-i,js} + \O (1)
\eeq
 as well as $| \hat{p}_{ijs} | \leq C (r_{j,s} )^{2/3}$. Therefore, for fixed $s, j$ we have, by applying Lemma \ref{lem:expectation-expand}
\begin{align}
& \sum_{i=1}^k f ( r_{j,s} - \hat{p}_{ijs} , r_{j, s} + \hat{p}_{ijs} )  = \sum_{i=1}^k r_{j,s} f ( 1 -  r_{j,s}^{-1} \hat{p}_{ijs} ,1  + r_{j, s}^{-1}  \hat{p}_{ijs} ) \notag  \\
= & r_{j,s} \sum_{i=1}^{k/2} f ( 1 -  r_{j,s}^{-1} \hat{p}_{ijs} ,1  + r_{j, s}^{-1}  \hat{p}_{ijs} ) + f ( 1 -  r_{j,s}^{-1} \hat{p}_{k+1-i,js} ,1  + r_{j, s}^{-1}  \hat{p}_{k+1-i,js} ) \notag \\
= &r_{j,s} k \mu +  \sum_{i=1}^{k/2} \mfa ( \hat{p}_{ij,s} + \hat{p}_{k+1-i,js} ) + \O ( k (r_{j,s} )^{1/3} ) \notag \\
= & r_{j, s } k \mu + \O ( k ( 2^j)^{1/2} )
\end{align}
where in the last line we applied \eqref{eqn:phase-2-a1} as well as the fact that $r_{j, s} \leq C (2^j)^{3/2}$. Therefore,
\beq
\left| (h_2 - h_1) \mu k - \sum_{j=1}^{N_2} \sum_{i=1}^k \sum_{s=1}^k f ( q_{ij,s}^{(2)} - q_{ij,s-1}^{(2)} )  \right|  \leq C \sum_{j=1}^{N_2} k^2 (2^j)^{1/2} \leq C k^{5/3} n^{1/3}
\eeq
as desired. \qed

The previous two lemmas immediately give the following.
\bep \label{prop:phase-2} There is are $C, c>0$ so that,
\beq
\pp\left[ \sum_{i=1}^k \log Z (q_{i,1}, q_{i,2} ) \leq \mu k (h_2 - h_1) - C k^{5/3} n^{1/3} \right] \leq C \e^{ - c k^2}.
\eeq
\eep

\subsection{Middle phase}

At the end of the previous phase, the $k$ paths intersect the line $\{ (x, y) : x +y = 2 h_2 \}$ on the $k$ points $q^{(3)}_{i,0}$ that have coordinates,
\beq
q^{(3)}_{i, 0} := \left( h_0 - p^{(3)}_i, h_0 + p^{(3)}_i \right)
\eeq
where
\beq
p^{(3)}_i := \left( 2 ( 2^{N_2} -1) + 2 \times 10^4 \right) \left( \frac{k+1}{2} - i \right) .
\eeq
Note that the separation between consecutive points along this line is of  order $n^{2/3} k^{-2/3}$. Set now $\ell^{(3)}_0 = h_2$ and 
\beq
\ell_j^{(3)} = \ell^{(3)}_{j-1} + \frac{ n - 2 h_2}{k}
\eeq
for $j=1, \dots, k$. We will demand that the $i$th curve passes through level $\ell_j$ at the point,
\beq
q^{(3)}_{ij} := \left( \ell^{(3}_j - p^{(3)}_i ,  \ell^{(3}_j + p^{(3)}_i \right) .
\eeq
We then bound
\begin{align}
\sum_{i=1}^k \log Z_{q_{i,2} ,q_{i,3} } & \geq \sum_{j=1}^k \sum_{i=1}^k \log Z ( q^{(3)}_{i,j-1} , q^{(3)}_{i,j} )  \notag \\
&\geq  \sum_{j=1}^k \sum_{i=1}^k \log \hat{Z}^{(3)} ( q^{(3)}_{i,j-1} , q^{(3)}_{i,j} )
\end{align}
where $\hat{Z}^{(3)} (x, y)$ is the partition function of polymer paths starting at $x$ and ending at $y$ that stay within $ c n^{2/3} k^{-2/3}$ of the straight line connecting $x$ to $y$, for a small enough $c>0$ so that the corridors of different paths are disjoint. 
\bel
There are $C, c>0$ so that,
\beq
\pp\left[  \sum_{j=1}^k \sum_{i=1}^k \left( \log \hat{Z}^{(3)} ( q^{(3)}_{i,j-1} , q^{(3)}_{i,j} )  - f( q^{(3)}_{i,j} - q^{(3)}_{i,j-1} ) \right) < C k^{5/3} n^{1/3} \right] \leq \e^{ - c k^2}
\eeq
and
\beq
 (h_3- h_2) k \mu = \sum_{j=1}^k \sum_{i=1}^k  f( q^{(3)}_{i,j} - q^{(3)}_{i,j-1} ) 
\eeq
\eel
\proof Let $Y_{ij} :=  f( q^{(3)}_{i,j} - q^{(3)}_{i,j-1} ) -\log \hat{Z}^{(3)} ( q^{(3)}_{i,j-1} , q^{(3)}_{i,j} )$. By the choice of the constraints, the random variables $\{ Y_{ij} \}_{ij}$ are all independent. The distance between $q^{(3)}_{i,j-1} $ and $q^{(3)}_{i,j}$ is of order $r:= (n/k)$ and the anti-diagonal displacement is $0$. The polymer paths are restricted to lie in a corridor of width of order $(n/k)^{2/3} \leq C r^{2/3}$ and so Proposition \ref{prop:constrained-lower-tail} is applicable. Therefore,
\beq
\pp\left[ Y_{ij} \geq  C r^{1/3} u \right] \leq \e^{ -c u}
\eeq
for $u \geq u_0$. So we may apply Proposition \ref{prop:conc} to the sum $\sum_{ij} Y_{ij}$. We have,
\beq
\nu \asymp k^2 r^{1/3} = k^{5/3} n^{1/3}, \qquad a_* \asymp r^{-1/3} = k^{1/3} n^{-1/3} ,
\eeq
and so,
\beq
\pp\left[ \sum_{ij} Y_{ij} \geq C k^{5/3} n^{1/3} \right] \leq \e^{ - c k^2}.
\eeq
This completes the estimate of the lemma. The second statement follows from the fact that $f( q^{(3)}_{i,j} - q^{(3)}_{i,j-1} ) = ( \ell^{(3)}_j - \ell^{(3)}_{j-1} ) \mu$. \qed

We therefore obtain,
\bep \label{prop:phase-3}
There are $C>0$ and $c>0$ such that
\beq
\pp \left[ \sum_{i=1}^k \log Z (q_{i,2} ,q_{i,3} ) \leq (h_3  -h_2) \mu k - C k^{5/3} n^{1/3} \right] \leq \e^{ - ck^2} .
\eeq
\eep

\subsection{Tail bound}

From all of the previous, we obtain the following.
\bep \label{prop:watermelon-bound} There is a $c_1 >0$ and $c, C>0$ so that for any $k = 2^N \leq c_1 n$ we have,
\beq
\pp\left[ \log Z_{(0, 0), (n, n)} \leq n \mu - C k^{2/3} n^{1/3} \right] \leq C \e^{  - ck^2} + C \e^{ - c k^{4/3} n^{2/3} \log(k)^{-1} }
\eeq
\eep
\proof We recall,
\beq
\log Z_{(0, 0), (n, n)} - n \mu \geq \frac{1}{k} \sum_{m=1}^5 \left( \sum_{i=1}^k \log Z ( q_{i,m-1} , q_{i,m} )  - (h_m - h_{m-1} ) n \mu \right). 
\eeq
The terms on the RHS with $m=1, 2, 3$ were bounded by Propositions \ref{prop:phase-1}, \ref{prop:phase-2} and \ref{prop:phase-3}, respectively. The terms with $m=4, 5$ are treated by the mirror opposite constructions of the phases $m=2, 1$, respectively. \qed

\section{Lower bound for upper tail} \label{sec:upper-tail}

 First, by convergence of $n^{-1/3} ( \log Z_{(0, 0), (n, n)} - \mu n )$ to a Tracy-Widom random variable \cite{BCF} and the unbounded support of this distribution we have that there is $c_1, \delta_1 >0$ so that
\beq
\pp\left[ \log Z_{(0, 0), (n, n)} \geq \mu n + c_1 n^{1/3} \right] \geq \delta_1.
\eeq
That the Tracy-Widom distribution has unbounded support can be deduced from the explicit form of its distribution function (cdf). In fact, more is known \cite[(25)]{baik-lower}:
\[F_2(x)=1-\frac{1}{32\pi x^{3/2}}\exp\Big(-\frac{4}{3}x^{3/2}\Big)(1+o(1)), \quad x\rightarrow \infty.\]
Here, $F_2(x)$ is the cdf of the Tracy-Widom GUE distribution. 

Fix some $k$. We have,
\beq
\log Z_{(0, 0), (n, n)} \geq \sum_{j=1}^k \log Z_{ \frac{j-1}{k} (n, n) , \frac{j}{k} (n, n)  }
\eeq
and by independence,
\beq
\pp\left[ \log Z_{(0, 0), (n, n) } \geq \mu n + u n^{1/3} \right] \geq \left( \pp\left[ \log Z_{(0, 0), (n/k, n/k ) } \geq \mu n/k + u n^{1/3}/k \right] \right)^k.
\eeq
We take $k^{2/3} = u/c_1$ so that with $r=n/k$ we have,
\beq
\pp\left[ \log Z_{(0, 0), (n/k, n/k ) } \geq \mu n/k + u n^{1/3}/k \right] = \pp\left[ \log Z_{(0, 0), (r, r) } \geq \mu r + c_1 r^{1/3} \right] \geq \delta_1
\eeq
for all $r$ large enough. The lower bound of \eqref{eqn:main-upper-tail} follows.
\appendix

\section{Concentration estimates}

\subsection{Sub-exponential random variables}

We have the following, Theorem 5.1(i) of \cite{janson}.
\bep \label{prop:exp-conc}
Let $W = \sum_{i=1}^n W_i$ where $W_i$ are independent exponential random variables $W_i \sim \mathrm{Exp} (a_i )$. Let,
\beq
\nu = \ee W = \sum_i \frac{1}{a_i}, \qquad a_* := \min_i a_i.
\eeq
Then for $\lambda \geq 1$ we have,
\beq
\pp\left[ W \geq \lambda \nu \right] \leq \lambda^{-1} \e^{ - a_* \nu ( \lambda - 1 - \log \lambda ) }.
\eeq
\eep

With this we prove the following.
\bep \label{prop:conc}
Let $\{Y_i \}_i$ be a collection of independent random variables such that for some $\theta_0$ and $\{ a_i \}_i$ we have,
\beq
\pp\left[ Y_i \geq \theta (a_i )^{-1} \right] \leq \e^{ - \theta }
\eeq
for all $i$ and $\theta \geq \theta_0$. Define,
\beq
\nu = \sum_{i=1}^n \frac{1}{ a_i}, \qquad a_* := \min_i a_i.
\eeq
Then there are $C>0$ and $c>0$ so that,
\beq
\pp\left[ \sum_i Y_i \geq C \nu \right] \leq  \e^{ - c a_* \nu }.
\eeq
\eep
\proof There is a coupling of $\{ Y_i \}_i$ to a family of mutually independent exponential random variables $X_i \sim \mathrm{Exp} (a_i)$ such that$^1$\let\thefootnote\relax\footnotetext{1. This can be constructed, e.g.,  by setting $Y_i = F^{-1} (U_i)$ where $F^{-1}$ is a generalized inverse of the CDF of $Y$ and $U_i$ are iid uniform $(0, 1)$ random variables. }
\beq
Y_i \leq X_i + a_i^{-1} \theta_0. 
\eeq
Taking $C \geq  \theta_0+3$ we then see that,
\beq
\pp\left[ \sum_i Y_i \geq C \nu \right] \leq \pp\left[ \sum_i X_i \geq 3 \nu \right] \leq \e^{ - c a_* \nu }
\eeq
where we applied Proposition \ref{prop:exp-conc} in the last inequality with $\lambda =3 $. \qed

\subsection{Sub-Gaussian random variables}

For a random variable $X$ we define the sub-Gaussian norm $ \| X \|_{\psi_2}$ by
\beq
\| X \|_{\psi_2} := \inf \left\{ K > 0 : \pp \left[ |X| >t \right] \leq 2 \e^{ - t^2 / K^2} , \forall t >0 \right\}
\eeq
For sub-Gaussian random variables we have the following, \cite[Theorem 2.6.3]{V}.
\bet \label{thm:v-thm} There is a $c>0$ so that the following holds.
Let $\{ X_i \}_{i=1}^N$ be mean-zero, independent sub-Gaussian random variables and let $K = \max_i \| X \|_{\psi_2}$. Let $a=(a_1, \dots, a_N) \in \rr^N$. Then,
\beq
\pp\left[ \left| \sum_{i=1}^N a_i X_i \right| > t \right] \leq 2 \exp \left(  - \frac{ c t^2}{ K^2 \| a\|_2^2} \right)
\eeq
\eet
As an application we have the following.
\bep \label{prop:sub-gauss}
Let $\{ G_i \}_{i=1}^N$ be a family of independent random variables such that there are $C_0, c_0>0$ so that
\beq
\pp\left[ |G_i | > t \right] \leq C_0 \e^{ - c_0 t^2}
\eeq
for $|t| > C_0$.  There are $C_1, c_1 >0$ depending only on $C_0, c_0 >0$ and not on $N$ so that for any $a = (a_1, \dots a_N) \in \rr^N$ we have,
\beq
\pp\left[ \left| \sum_{i=1}^N a_i G_i \right| > C_1 \| a \|_1 + t \right] \leq C_1 \exp \left( - \frac{ c_1 t^2}{ \| a \|_2^2} \right)
\eeq
\eep
\proof Define $X_i := G_i - \ee[ G_i]$. Since $| \ee[ G_i ] | \leq C$ for all $i$, we see that there is a $c>0$ depending only on $c_0, C_0 >0$ so that
\beq
\pp\left[ | X_i | > t \right] \leq 2 \e^{ - c t^2}
\eeq
and so $K:= \max_i \| X_i \|_{\psi_2}$ is bounded by a constant depending only on $c_0, C_0 >0$. From Theorem \ref{thm:v-thm} we see that,
\beq
\pp\left[ \left| \sum_{i=1}^N a_i X_i \right| > t \right] \leq 2 \exp \left(  - \frac{ c_1 t^2}{  \| a\|_2^2} \right)
\eeq
for some $c_1 >0$. On the other hand we have that
\beq
\left| \ee\left[ \sum_{i=1}^N a_i G_i \right] \right| \leq C_1 \|a \|_1
\eeq
for some $C_1 >0$ 
and so the claim follows. \qed
\section{FKG inequality}
In this section we prove a form of positive association (``the Harris-FKG inequality") for polymer partition functions. For $1 \leq i \leq N$ and $1 \leq j \leq M_i$, let $X_{ij}$ be a random variable of the form
\beq
X_{ij} = \int_{ a_{ij} < s_{m_{ij}} < \dots s_{n_{ij}} < b_{ij} } \e^{ \sum_{k=m_{ij}}^{n_{ij}} B_k (s_k ) - B_{k} (s_{k-1 } ) } \prod_{k=m_{ij}}^{n_{ij} -1} \1_{ \{ s_k \in I_{ijk} \} } \d s_{m_{ij}} \dots \d s_{n_{ij}-1}
\eeq
where the $I_{ijk}$ are some intervals. Then, let
\beq \label{eqn:Zi-def}
Z_i = \sum_{j=1}^{M_i} a_{ij} X_{ij}
\eeq
where $a_{ij} >0$. 
\bep \label{prop:fkg}
Let the $Z_i$ be as above. Let $A_i \in \rr$. Then,
\beq
\pp\left[ \bigcap_{i=1}^N \{ Z_i \leq A_i \} \right] \geq \prod_{i=1}^N \pp \left[ Z_i < A_i \right] 
\eeq
\eep
\proof Let $L \geq \max_{ij} | n_{ij}| + \max_{ij} | m_{ij}  | + \max_{ij} |a_{ij}|  + \max_{ij} | b_{ij} |$. 
Let $\{ Y_{kl} \}_{(k, l) \in \zz^2}$ be a family of iid $\pm 1 $ random variables. Consider for every $n$ the functions $\hat{B}^{(n)}_i : [-L, L] \to \rr$ defined by
\beq \label{eqn:hatB-def}
\hat{B}^{(n)}_i (t) := n^{1/2} \int_{-L}^t \sum_{j} Y_{ji} \1_{ \{ s \in (j/n, (j+1)/n ) \} } \d s .
\eeq
Then each $B^{(n)}_i (t)$ converges in the space $C([-L, L] )$ equipped with the topology induced by the uniform norm $ \| \cdot \|_\infty$ to Brownian motions $W_i (t)$ on $[-L, L]$ with $W_i (-L) = 0$. Viewing the random variables $X_{ij}$ as functions $X_{ij} (\cdot ) : C([-L, L] )^{2L+1} \to \rr$, we see that they are continuous with respect to the norm $\| \cdot \|_\infty$. This implies the joint convergence of the collection $\{ Z_i ( B^{(n)} ) \}_i$ to $\{ Z_i (W ) \}_i$. However, this latter family has the same distribution as the original $Z_i$ specified in \eqref{eqn:Zi-def} in terms of the original Brownian motions $B_i (t)$. By the Portmanteau theorem, 
\beq
\limsup_{n \to \infty} \pp\left[ \bigcap_i \{ Z_i (B^{(n)} ) \leq A_i \} \right] \leq \pp\left[ \bigcap_i \{ Z_i \leq A_i \} \right]
\eeq
and
\beq
\liminf_{n \to \infty} \prod_i \pp\left[ \{ Z_i ( B^{(n)} ) < A_i \} \right] \geq \prod_i \pp\left[ \{ Z_i < A_i \} \right] .
\eeq
On the other hand, from the representation \eqref{eqn:hatB-def}, we see that any increment $\hat{B}_i^{(n)} (t) - \hat{B}_i^{(n)} (s) $ is increasing under changing any $Y_{ij}$ from $-1$ to $+1$. Therefore, by positive association (the Harris-FKG inequality for independent random variables, see, e.g., \cite[Chaper II.2]{Lig}),
\beq
\prod_i \pp\left[ \{ Z_i ( B^{(n)} ) < A_i \} \right] \leq \pp\left[ \bigcap_i \{ Z_i (B^{(n)} ) \leq A_i \} \right]
\eeq
for every $n$. The claim follows. \qed

\section{Miscellaneous proofs}

In this section it will be useful to introduce the notation,
\beq
B_k (t, s) := B_k (t) - B_k (s)
\eeq
for the Brownian increments.

\subsection{Proof of \eqref{eqn:fd}} \label{sec:fd}

Recall that $T_{jk}^{(+)}$ is the set of polymer paths intersectings the lines $\ell_i := \{ (z_i, z_i ) + (-m , m) : |m| \leq b_{j-1} n^{2/3} \}$ for $i=1, 2$ and passing above the line $\{ (z_0, z_0 ) + (-m, m) : |m| \leq b_j n^{2/3} \}$. Write $a = b_{j-1} n^{2/3}$ and $b = b_j n^{2/3}$. The event that the polymer path intersects the line $\ell_i$ can be written as the disjoint union of the sets (up to some sets of Lebesgue measure $0$ which do not contribute to the partition function)
\beq
\left( \bigsqcup_{|m| \leq a } \{ s_{z_i + m } > z_i - m , s_{z_i + m -1 } < z_i - m \} \right) \bigsqcup \left( \bigsqcup_{-a \leq m \leq a-1 } \{ z_i -m -1 < s_{z_i + m } <z_i - m \}  \right) 
\eeq
and the event that the path passes above the line $\{ (z_0, z_0 ) + (-m, m) : |m| \leq b \}$ can be written as the event $\{ s_{z_0 + b } < z_0 - b \} =: \A$. Therefore,
\begin{align}
 & Z_{(0, 0), (n, n) } [ T_{jk}^{(+)} ] \notag \\
 =& \sum_{|i|, |j| \leq a } Z_{(0, 0), (n, n)} [ s_{z_1+i} > z_1 - i , s_{z_1+i-1} < z_1 - i , s_{z_2+j} > z_2 - j , s_{z_2+j-1} < z_2 - j ,  \A] \label{eqn:fd-1} \\
 +& \sum_{-a \leq i , j \leq a-1 } Z_{(0, 0), (n, n)} [ z_1 - i - 1 < s_{z_1 + i } < z_1 - i , z_2 - j - 1 < s_{z_2 + j } < z_2 - j, \A] \label{eqn:fd-2} \\
  +& \sum_{ |i| \leq a , -a \leq j \leq a-1 } Z_{(0, 0), (n, n)} [ s_{z_1+i} > z_1 - i , s_{z_1+i-1} < z_1 - i ,z_2 - j - 1 < s_{z_2 + j } < z_2 - j, \A]  
  \label{eqn:fd-3} \\
  +& \sum_{ -a \leq i \leq a-1, |j| \leq a } Z_{(0, 0), (n, n)} [ z_1 - i - 1 < s_{z_1 + i } < z_1 - i , s_{z_2+j} > z_2 - j , s_{z_2+j-1} < z_2 - j ,  \A ]
\label{eqn:fd-4}
\end{align}
The terms on the first, second, third and fourth lines above will be seen to give the first, second, third and fourth terms in \eqref{eqn:fd}, respectively. 
For the terms on the line \eqref{eqn:fd-1}, the set
\beq
\{ s_{z_1+i} > z_1 - i , s_{z_1+i-1} < z_1 - i , s_{z_2+j} > z_2 - j , s_{z_2+j-1} < z_2 - j  , s_{z_0+b} < z_0 - b \} 
\eeq
is empty unless $z_1 - i < z_0 - b$ and $z_2 + j > z_0 + b$. Since $a \leq b$ and $|i| \leq a$ the first inequality implies $z_0 +b > z_1 +i$. For the non-zero terms, group the terms in the integrand as,
\begin{align}
 & \left(  \1_{ \{ s_{z_1+i -1} < z_1 - i \} } \e^{ \sum_{k=0}^{z_1+i-1} B_k (s_k, s_{k-1} ) + B_{z_1+i} (z_1 - i, s_{z_1+i-1} ) } \prod_{k=0}^{z_1+i-1} \d s_{k} \right) \notag \\
\times & \bigg\{  \1_{ \{ s_{z_0+b} < z_0 - b \} } \1_{ \{ s_{z_1 + i} > z_1 - i \} } \1_{ \{ s_{z_2+j-1} < z_2 - j \} } \notag \\
\times & \e^{ B_{z_1+i} (s_{z_1+i},z_1 - i ) + B_{z_2+j} ( z_2 - j , s_{z_2 + j -1 } ) + \sum_{k=z_1 + i + 1}^{ z_2 +j -1 } B_k (s_k , s_{k-1} ) } \prod_{k=z_1  +i}^{z_2 +j -1 } \d s_k  \bigg\} \notag \\
& \times \left( \1_{ \{  s_{z_2+j } > z_2 - j \} } \e^{ B_{z_2 + j} ( s_{z_2+j} , z_2 - j ) + \sum_{k=z_2+j+1}^n B_k (s_k, s_{k-1} ) } \prod_{k=z_2+j}^{n-1} \d s_k \right)
\end{align}
so that,
\begin{align}
 & \sum_{|i|, |j| \leq a } Z_{(0, 0), (n, n)} [ s_{z_1+i} > z_1 - i , s_{z_1+i-1} < z_1 - i , s_{z_2+j} > z_2 - j , s_{z_2+j-1} < z_2 - j ,  \A]  \notag \\
= & \sum_{ |i|, |j| \leq a } Z_{(0, 0), (z_1 - i, z_1 + i ) } Z_{(z_1 - i, z_1+i), (z_2 -j, z_2 + j ) } [ \A] Z_{ (z_2 - j , z_2 + j ) , (n, n) } 
\end{align}
where the extra terms on the RHS that correspond to terms we argued were zero on the LHS are automatically $0$ because $Z_{(z_1 - i, z_1+i), (z_2 -j, z_2 + j ) } [ \A] = 0$ if either $z_1 - i \geq z_0 - b$ or $z_2 +j \leq z_0 + b$. 

For the terms on the line \eqref{eqn:fd-2}, the set
\beq
\{ z_1 - i - 1 < s_{z_1 + i } < z_1 - i , z_2 - j - 1 < s_{z_2 + j } < z_2 - j,s_{z_0+b} < z_0 - b \}
\eeq
is empty unless $z_1 - i \leq z_0 - b$ and $z_2+j > z_0 + b$. The first inequality implies $z_0 + b > z_1 + i$ since $b \geq a \geq i +1 $. For the non-zero terms on the line \eqref{eqn:fd-2} we group the integrand as,
\begin{align}
 &\bigg\{  \left( \1_{ \{  s_{z_1 + i-1} < s_{z_1 + i } \} } \e^{ \sum_{k=0}^{z_1 + i } B_k (s_k , s_{k-1} ) } \prod_{k=0}^{z_1 +i-1} \d s_k \right) \notag \\
 \times & \left( \1_{ \{s_{z_2 +j+1 } > s_{z_2 + j} \} } \e^{ \sum_{k=z_2+j+1}^{n} B_k (s_k , s_{k-1} ) } \prod_{k=z_2+j+1}^{n-1}  \d s_k \right) \notag \\
 \times & \left( \1_{ \{ s_{z_1 +i +1} > s_{z_1 + i } \} } \1_{ \{ s_{z_0+b} < z_0 + b  \} } \1_{ \{ s_{z_2 + j -1 } < s_{z_2 + j } \} } \e^{ \sum_{k=z_1+i+1}^{z_2+j}  B_k (s_k , s_{k-1} ) } \prod_{k=z_1+i+1}^{z_2+j-1} \d s_k \right)   \bigg\} \notag \\
& \1_{  \{ z_1 - i - 1 < s_{z_1 +i} < z_1 - i \}  } \1_{  \{ z_2 - j - 1 < s_{z_2+j} < z_2 - j  \} } \d s_{z_1+i} \d s_{z_2 +j } 
\end{align}
Therefore, the sum on the second line \eqref{eqn:fd-2} equals,
\begin{align}
\sum_{ a \leq i , j \leq a -1 } \int & \bigg\{  Z_{(0, 0), (s_{z_1+i}, z_1 + i ) }  Z_{(s_{z_1+i}, z_1 + i + 1 ) , (s_{z_2 + j } , z_2+j ) } [ \A] Z_{(s_{z_2+j} , z_2 + j + 1 ), (n, n) }  \notag \\
 & \times \1_{  \{ z_1 - i - 1 < s_{z_1 +i} < z_1 - i \}  } \1_{  \{ z_2 - j - 1 < s_{z_2+j} < z_2 - j  \} } \bigg\} \d s_{z_1+i} \d s_{z_2 +j } ,
\end{align}
where again, terms that were zero on the line \eqref{eqn:fd-2} are also zero above. The remaining lines \eqref{eqn:fd-3} and \eqref{eqn:fd-4} are handled via highly similar arguments which are omitted. \qed

\subsection{Proof of \eqref{eqn:lb-a2}} \label{sec:sd}

Let us label $L_i := \{ (z_i, z_i ) - (m, -m ) : |m| \leq a \}$. The event that the polymer path intersects the line $L_i$ can be written as the disjoint union of the events (up to sets of Lebesgue measure $0$ which do not contribute to the partition function),
\begin{align}
& \left( \bigsqcup_{ |m| \leq a } \{ s_{z_i + m } > z_i - m , s_{z_i-m -1 } < z_i - m \} \right) \bigsqcup \left( \bigsqcup_{- a \leq m \leq a-1 } \{ z_i - m -1 < s_{z_i+m} < z_i -m \} \right) \notag \\
& =:\left( \bigsqcup_{ |m| \leq a } \A^{(0)}_{m,i} \right) \bigsqcup \left( \bigsqcup_{- a \leq m \leq a-1 } \A^{(1)}_{m, i} \right) =: \A^{(0)}_{i} \sqcup \A^{(1)}_{i} .
\end{align} 
Therefore, 
\beq
Z_{(0, 0), (n, n) }[ \A] = \sum_{ \sigma \in \{ 0, 1 \}^{k-1} } Z_{(0, 0), (n, n) } \left[ \bigcap_{i=1}^{k-1} \A^{( \sigma_i ) }_i \right]
\eeq
A single term in the summand on the RHS can be written,
\beq
Z_{(0, 0), (n, n) } \left[ \bigcap_{i=1}^{k-1} \A^{( \sigma_i ) }_i \right] = \sum_{ \substack{ m_1, \dots m_k  \\  -a \leq m_i \leq a - \sigma_i, \forall i }} Z_{(0, 0), (n, n)} \left[ \bigcap_{i=1}^{k-1} \A_{m_i, i}^{(\sigma_i ) } \right]
\eeq
Each term on the RHS is non-zero only if  $z_i + m_i \geq z_{i-1} + m_{i-1} + \sigma_{i-1}$ for $i=1, k$ with the convention $z_0+m_0 + \sigma_0 = 0$ and $z_k + m_k + \sigma_k = n$. For such terms we have the decomposition,
\begin{align}
& Z_{(0, 0), (n, n)} \left[ \bigcap_{i=1}^{k-1} \A_{m_i, i}^{(\sigma_i ) } \right] \\
=& \int Z_{(0, 0), (x_1, y_1) } Z_{( x_{k-1}, y_{k-1} + \sigma_{k-1} ), (n, n) } \prod_{i=1}^{k-2} Z_{ ( x_i, y_i + \sigma_i ) , (x_{i+1}, y_{i+1} ) } \prod_{i=1}^{k-1} \d \xi_i (x_i, y_i )
\end{align}
where $\xi_i$ is a delta function at $(z_i - m_i, z_i + m_i)$ if $\sigma_i = 0$ and $1$d Lebesgue measure on the interval $\{ (z_i - m_i, z_i +m_i ) - (s, 0) : 0 < s < 1 \}$ if $\sigma_i =1$. Note that the above identity extends to the excluded case where the LHS is $0$ as so is the RHS by inspection (recall our convention $Z_{p, q} = 0$ if the coordinate-wise ordering $p \leq q$ does not hold). \qed

\subsection{Proof of \eqref{eqn:lb-a1}} \label{sec:misc-1}
 On the complement of the event of \eqref{eqn:deviation-event} we have by \eqref{eqn:deviation-main} that for $ u + (j-1) n^{-10} \leq w \leq u + j n^{-10}$ that
 \beq
 Z_{(s, m), (w, p) } \leq \e Z_{(s, m), (u+j n^{-10} , p ) } .
 \eeq
 Similarly,
 \beq
 Z_{(w, p+1), (t, q) } \leq \e Z_{(u+(j-1)n^{-10} , p+1 ), (t, q) } .
 \eeq
 The claim then follows. \qed

 \subsection{Proof of \eqref{eqn:tv-a1}} \label{sec:misc-2}
 
 In order to prove this inequality we use the representation \eqref{eqn:fd}. There are several components. First, we clearly have,
 \begin{align}
  & \int \int \tilde{Z}_{(0, 0), p} \tilde{Z}_{p, q} [  \A ]  \tilde{Z}_{q, (n, n)}\d \mu_1 (p) \d \mu_1 (q) \\
 \leq & \left( \int \tilde{Z}_{(0, 0), p } \d \mu_1 (p ) \right) \left( \int \tilde{Z}_{p, q} [ \A ] \d \mu_1 (p ) \d \mu_2 (q ) \right) \left( \int \tilde{Z}_{q, (n, n) } \d \mu_2 (q ) \right)
 \end{align}
 since the measures $\mu_i$ are simply sums of delta functions. We now consider the second term of \eqref{eqn:fd}, the integral against $\d \nu_1 \d \nu_2$.  The event,
 \beq \label{eqn:deviation-app}
\B := \bigcap_{k=-10n}^{10n} \left \{ \sup_{|s_1|, |s_2| \leq 10 , i=1, 2} |B_{ z_i -k } (z_i -k +s_1 ) - B_{z_i -k} (z_i - k + s_2 ) | \leq \delta  b^2 n^{1/3} \right\}
 \eeq
 holds with probability at least $ 1- C n \e^{ - c \delta^2 b^4 n^{2/3} }$. On this event we see from \eqref{eqn:deviation-main} that for $s \in (0, 1)$ and $|m| \leq n$,
 \beq \label{eqn:misc-2-a1}
 Z_{(0, 0), (z_1 - m -s, z_1 +m ) } \leq \e^{  \delta b^2 n^{1/3} } Z_{(0, 0), (z_i - m , z_i + m ) }
 \eeq
 as well as
 \beq \label{eqn:misc-2-a2}
 Z_{ (z_2 -m - s , z_2 +m + 1 ), (n, n) } \leq \e^{ \delta b^2 n^{1/3} } Z_{ (z_i - m -1, z_i + m + 1 ), (n, n) }.
 \eeq
 For $|i|, |j| \leq b_{j-1} n^{2/3}$ and $t_1, t_2 \in (0, 1)$ the term,
 \beq
 Z_{(z_1 - i - t_1 , z_1 +i +1 ) , (z_2 - j - t_2 , z_2 + j ) } [ \A ] 
 \eeq
 is $0$ unless $z_1 - i \leq z_0 - b_j n^{2/3}$ (since $\A = \{ s_{z_0 + b_j n^{2/3} } < z_0 - b_j n^{2/3} \}$), which implies $z_0 + b_j n^{2/3} > z_1 + i + 2$ (since $b_j > b_{j-1} + 10$ and $i \leq b_{j-1} n^{2/3}$). This term is also $0$ unless $z_1 + j > z_0 + b_j n^{2/3}$. We then have the representation,
 \begin{align}
  & Z_{(z_1 - i - t_1 , z_1 +i +1 ) , (z_2 - j - t_2 , z_2 + j ) } [ \A ] \notag \\
  = & \int_{z_1 - i -t_1 < s_{z_1 +i + 1 } < s_{z_2+j-1} < z_2 - j - t_2 }  \bigg\{ Z_{(s_{z_1 + i + 1  } , z_1 + i + 2 ) , ( s_{z_2 +j-1 } , z_2 + j -1 )} [ \A_1] \1_{ \{  s_{z_2+j-1} \in \A_2  \} } \notag \\
  \times &  \e^{ B_{z_1 + i + 1 } ( s_{z_1+i+1} , z_1 - i - t_1 ) + B_{z_2+j} ( z_2 - j - t_2 , s_{z_2+j-1} ) } \bigg\} \d s_{z_1 + 1 + i } \d s_{z_2 + j -1 } 
 \end{align}
 where if $z_2 + j -1 = z_0 + b_j n^{2/3}$ we have $\A_1 $ is the whole polymer space and $\A_2 = (-\infty, z_0 - b_j n^{2/3} )$ and otherwise $\A_1 = \A$ and $\A_2 = \rr$. From this representation we conclude that on the event $\B$ that,
 \beq \label{eqn:misc-2-a3}
  Z_{(z_1 - i - t_1 , z_1 +i +1 ) , (z_2 - j - t_2 , z_2 + j ) } [ \A ] \leq \e^{ 2 \delta b^2 n^{1/3} } Z_{ (z_1 - i - 1 , z_1 + i + 1 ) , (z_2 - j , z_2 + j ) } [ \A]
 \eeq
 in a similar manner to the proof of Lemma \ref{lem:deviation}. By combining \eqref{eqn:misc-2-a1}, \eqref{eqn:misc-2-a2} and \eqref{eqn:misc-2-a3} we see that the term on the second and third lines of \eqref{eqn:fd} is bounded by
 \beq
 C \e^{ 4 \delta b^2 n^{2/3} } \int \int \tilde{Z}_{(0, 0), p} \tilde{Z}_{p, q} [ \A] \tilde{Z}_{q, (n, n) } \d \mu_1 (p) \d \mu_2 (q) .
 \eeq
 The other terms of \eqref{eqn:fd} can be estimated in a similar fashion. The claim follows after sending $\delta = c \delta$, for some $c>0$. \qed


\end{document}